\documentclass[AMA,STIX1COL]{WileyNJD-v2}

\articletype{Article Type}%

\usepackage{tikz}
\usetikzlibrary{arrows, arrows.meta} 
\usetikzlibrary{decorations.pathreplacing}
\usepackage{adjustbox}
\usetikzlibrary{decorations.pathmorphing,patterns}
\usetikzlibrary{calc,patterns,decorations.markings}
\usetikzlibrary{positioning}
\usepackage{float}
\usepackage{multicol}
\usepackage{amsmath, amsthm, amssymb}
\usepackage{titlesec}
\usepackage{paralist}
\usepackage{wrapfig}
\usepackage{setspace}
\titlespacing{\section}{0pt}{\parskip}{-\parskip}
\titlespacing{\subsection}{0pt}{\parskip}{-\parskip}
\titlespacing{\subsubsection}{0pt}{\parskip}{-\parskip}
\usepackage{subcaption}
\usepackage{blindtext}
\usepackage{paralist}
\usepackage{stackengine}

\newenvironment{rcases}
  {\left.\begin{aligned}}
  {\end{aligned}\right\rbrace}
\numberwithin{equation}{section}
\usepackage[section]{placeins}

\graphicspath{{Figs/}}
\newcommand{\ie}{\textit{i}.\textit{e}., }

\def\etc.{etc,.\spacefactor=\the\sfcode`\c}

\received{26 April 2016}
\revised{6 June 2016}
\accepted{6 June 2016}

\begin{document}

\title{Robust and accurate central algorithms for Multi-Component mixture equations with Stiffened gas EOS}

\author[1,4]{Ramesh Kolluru}

\author[2]{S V Raghurama Rao}

\author[3]{G N Sekhar}

\authormark{Ramesh Kolluru \textsc{et al}}

\address[1]{\orgdiv{Department of Aerospace Engineering}, \orgname{Indian Institute of Science}, \orgaddress{\state{Bangalore,Karnataka}, \country{India}}}

\address[2]{\orgdiv{Department of Aerospace Engineering}, \orgname{Indian Institute of Science}, \orgaddress{\state{Bangalore,Karnataka}, \country{India}}}

\address[3]{\orgdiv{Department of Mathematics}, \orgname{BMS College of Engineering}, \orgaddress{\state{Bangalore,Karnataka}, \country{India}}}

\address[4]{\orgdiv{Department of Mechanical Engineering}, \orgname{BMS College of Engineering}, \orgaddress{\state{Bangalore,Karnataka}, \country{India}}}
\corres{*Ramesh Kolluru  \email{kollurur@iisc.ac.in}}

\presentaddress{Post Doctoral Fellow, Department of Aerospace Engineering, Indian Institute of Science, Bangalore, Karnataka, India}

\abstract[Summary]{Simple and robust algorithms are developed for compressible Euler equations with stiffened gas equation of state (EOS), representing gaseous mixtures in thermal equilibrium and without chemical reactions. These algorithms use fully conservative approach in finite volume frame work for approximating the governing equations. Also these algorithms used central schemes with controlled numerical diffusion for this purpose. Both Mass fraction ($Y$) and $\gamma$ based models are used with RICCA and MOVERS+ algorithms to resolve the basic features of the flow fields. These numerical schemes are tested thoroughly for pressure oscillations and preservation of the positivity of mass fraction at least in the first order numerical methods. Several test cases in both 1D and 2D are presented to demonstrate the robustness and accuracy of the numerical schemes.
}


\jnlcitation{\cname{%
\author{Ramesh Kolluru.}, 
\author{S V Raghurama Rao}, and
\author{G N Sekhar}, 
} (\cyear{2019}), 
\ctitle{Simple,Accurate and Robust Algorithms for Multi-Component mixture equations with Stiffened gas EOS}, \cjournal{}, \cvol{}.}

\maketitle

\footnotetext{\textbf{Abbreviations:} MOVERS, Method of Optimum Viscosity for Enhanced Resolution of Shocks; RICCA Riemann Invariants based Contact-discontinuity Capturing Algorithm; }

\section{Introduction}
Atmospheric air is a mixture of gases which are compressible in nature. Each of the components in the mixture have different physical and thermodynamical properties, and very often in modelling the flow, air is assumed to be a single component gas with constant properties. 
There are many applications where due consideration should be given to each of the components in the mixture such as gasoline and air mixture entering the combustion chamber and combustion products exhausting from the engines. There are instances where liquids and gases exist together like bubbles moving in the liquid, spray of paint facilitated through a nozzle.  In all the situations mentioned above fluids exist as a mixture or as different components separated by interfaces. 

Many times, the contribution of the individual components are negligible or the variation in the properties of the components do not contribute significantly to the flow field and hence they can be neglected.  If the properties of the components vary at large, then the individual effect of the components are resolved or their combined effect on the mixture has to be studied. In these cases, the classic model of single component compressible fluid may not be appropriate. 

Broadly the flow of these fluids can be classified into two categories: a) pure interface problems, and b) multicomponent flows. For pure interface problems, the thermodynamic properties of the fluids change only across the interface whereas in multicomponent flows the properties vary throughout the flow field. In pure interface problems, apart from solving for the dynamics of each component, the interface is also tracked by a specific method like level set method. In multicomponent flows the modelling is done without tracking any interface. Nature of multicomponent flows can vary from low subsonic flows to hypersonic reacting flows. Low subsonic flows often are coupled with combustion related phenomena and therefore are not easily amenable to numerical modelling. Modelling of supersonic and hypersonic flows can take advantage of the sophisticated numerical methods developed for hyperbolic systems in the past few decades, though treating supersonic and hypersonic combustion problems are non-trivial. 
Some of the important contributions in modelling multicomponent flows are due to
\cite{karni1,karni2,karni3,Chargy,RA1,RA2,RA3,RA4,Larrouturou1,Larrouturou2,shyue1,shyue2,shyue3}.  

Fernandez \etal \cite{Fernandez}, aimed at constructing an efficient conservative numerical scheme for computation of multi-species flows. The governing equations are Euler equations and additional equations for the species with different molecular weights and specific heats are considered. Approximate Riemann solver of Roe has been used and modifications for evaluation of $\gamma$ in the Roe matrix have been suggested. Donor cell approximation method for species equations is modelled and compared with the modifications of the Roe scheme. They conclude that the modified Roe scheme performs better than the donor cell approximation method. 

Larrouturou \etal \cite{Larrouturou1}, have reviewed various numerical methods for multicomponent perfect and real gas models. They have suggested modifications for Osher, Steger-Warming, van Leer and Roe schemes for the application to multicomponent perfect and real gases. They show clearly that, for multicomponent flows Roe's conditions (consistency, conservation and hyperbolicity) get satisfied only when $\gamma(U)$ is constant and hence the extension of Roe scheme for multicomponent mixture flows is not an easy task.

Karni \cite{karni1} has carried out modelling of multicomponent fluids using Euler equations with an additional equation for the species. Both conservative and non-conservative form of the equations are considered and primitive form of the equations are recommended to avoid the pressure oscillations occurring near the material interface. Four different models of the governing equations with variable $\gamma$ and a level set method based on distance function are used in both conservative and primitive form. Roe linearisation method is used in numerical simulation and compared with second order upwind methods.  It is concluded that any fully coupled conservative based numerical scheme leads to pressure oscillations and non-preserving of positivity of mass fractions. She recommended the use of primitive variable based approach in order to avoid the pressure oscillations. The use of non-conservative form, however, leads to conservative errors and incorrect shock positions.

Abgrall \cite{RA1}, has used a quasi-conservative approach for the calculations for multicomponent cases and proved that an additional evolution equation of $\gamma$ and in particular of the form $\frac{1}{\left( \gamma - 1 \right) }$, is suggested to preserve mass fraction positivity and to avoid pressure oscillations. 

Abgrall and Karni \cite{karni3} have reviewed numerical algorithms commonly used in the simulations of multicomponent compressible fluid flow. They conclude that if separate equations for individual species are solved along with the mixture equations, then the numerical scheme developed preserves pressure equilibrium and mass fraction positivity. 

Keh-Ming Shyue in \cite{shyue1,shyue2,shyue3,shyue4}, has utilized Abgrall's model \cite{RA1} for compressible multicomponent flow problems using stiffened gas EOS, van der Waals EOS, Mie-Gruneisen EOS, Tait EOS.

Overall the basic issues in extending the single fluid conservative numerical schemes to multicomponent flows are
\begin{enumerate}
\item preserving the positivity of mass fraction,
\item avoiding pressure oscillations even in the first order numerical scheme,
\item difficulties in extension to more than 2 components.
\end{enumerate}
In this work novel and accurate central solvers MOVERS-n, MOVERS-1, developed by \cite{Jaisankar_SVRRao} along with MOVERS+, RICCA as explained in \cite{rkolluru} are applied to multicomponent flows to address some of the above issues.
\section{Governing equations for mixture with two components or species}
From the literature it is observed that there are many different ways in which the governing equations can be formulated. A simple case of non-reacting mixture equations with two components and without diffusion is considered in the present work. Two models based on mass fraction $\emph{Y}$ and $\gamma$ are chosen to test the algorithms in conservative cell centered finite volume frame work. In the following sections the governing equations for these models and basic algorithms used to discretise them are discussed briefly.
\subsection{Mass fraction based model}
Consider the mixture of gasses consisting of two species with following mixture properties:  pressure p, density $\rho$, velocity u and temperature $T$. 
The mixture pressure is given by Dalton's Law $p = p_1 + p_2$, mixture density $\rho = \rho_1 + \rho_2$, the mass fraction of the species $Y_k = \frac{\rho_k}{\rho},~ k = 1,2$.
Specific heat at constant pressure and constant volume of individual species, $c_{p_k}, c_{v_k},~ k = 1,2$ are considered to be constant, and the ratio of specific heats of individual species is given by  $\gamma_k 
=\frac{{c_p}_k}{{c_v}_k},~ k = 1,2$. The governing equations for the mixture in conservation form are given by (\ref{mixture_governing_equations})
\begin{align}
\begin{rcases}
\label{mixture_governing_equations}
\frac{\partial \rho}{\partial t} + \frac{\partial\left( \rho u\right)}{\partial x} &= 0, \\
\frac{\partial \rho u}{\partial t} + \frac{\partial \left( \rho u^2 + p \right)}{\partial x}&=0,\\
\frac{\partial \rho E_t}{\partial t} + \frac{\partial \left[ \left( \rho E_t + p \right) u \right]}{\partial x}&=0,\\
\frac{\partial \left(\rho Y_k\right)}{\partial t} + \frac{\partial \left( \rho Y_k u\right)}{\partial x}&=0,\qquad k = 1,2 \\
Y_1 + Y_2 &= 1
\end{rcases}
\end{align}
The value of ratio of specific heats,  $\gamma$, for the mixture, is defined as
$\gamma = \frac{{c_p}_{mixture}}{{c_v}_{mixture}} = \frac{\sum_k Y_k\gamma_k {c_v}_k}{\sum_k Y_k {c_v}_k} $ and the mixture pressure is given by $p = (\gamma -1)(\rho E_t - \frac{\rho u^2}{2})$. The equation of state for each individual component can be described by a function $p = p(\rho, e)$.
These governing equations are represented in compact notation as in the first equation of (\ref{cfd_notation}) where  enthalpy of the mixture is given by $H = E_t + \frac{p}{\rho}$. The ratio of specific heats for the mixture, $\gamma$ is a function of the conserved variable vector $U$, as $\gamma = \frac{U_4 \gamma_1 {c_v}_1 + (U_1 - U_4) \gamma_2 {c_v}_2}{U_4 {c_v}_1 + (U_1 - U_4) {c_v}_2} = \gamma\left(U_1,U_4\right)$. This property of $\gamma$ for the mixture plays a role in determining the hyperbolicity of the governing equations. It can also be observed that (\ref{cfd_notation}) is extension of Euler equations with an additional equations for the mass fraction of individual component gases. Hence if this set of governing equations satisfies the hyperbolicity principle then all the algorithms which are designed for Euler equations can be in principle extended to multicomponent fluids.
\begin{align}
\label{cfd_notation}
\frac{\partial U}{\partial t} + \frac{\partial F\left(U \right)}{\partial x} &=0
\end{align}
\begin{align}
U = \begin{bmatrix}
         U_1 \\ 
         U_2\\ 
         U_3 \\ 
         U_4
        \end{bmatrix} =  \begin{bmatrix}
         \rho \\ 
         \rho u\\ 
         \rho E \\ 
         \rho Y_k
        \end{bmatrix} 
        ;
        F(U) =  \begin{bmatrix}
         \rho u \\ 
         \rho u^2 + p \\ 
          (\rho E + p)u \\ 
          \rho Y_k u 
        \end{bmatrix} &= \begin{bmatrix}
         U_2 \\ 
         \frac{\left(3-\gamma\right)}{2}\frac{U_2^2}{U_1} + \left(\gamma - 1\right)U_3\\ 
         \frac{\left(3-\gamma \right)}{2} \frac{U_3 U_2}{U_1} - \frac{\left(\gamma -1\right)}{2} \frac{U_2^3}{U_1^2} \\
         \frac{U_4 U_1}{U_2}
        \end{bmatrix}
 \end{align}
\subsection {Hyperbolicity and eigenstructure for the mixture model}
To demonstrate the hyperbolicity of equations (\ref{cfd_notation}), it is required to evaluate the flux Jacobian matrix, it's eigenvalues and corresponding eigenvectors. The flux Jacobian matrix of the governing equations is given by (\ref{Mixture_Flux_Jacobian})
\begin{align}
\label{Mixture_Flux_Jacobian}
A(U)  &=\begin{bmatrix}
         0 & 1 & 0& 0 \\ 
          \frac{\left(\gamma-3\right)}{2}\frac{U_2^2}{U_1^2} + B & \left(3 - \gamma \right)\frac{U_2}{U_1}&\left(\gamma -1\right) &B'\\ 
          \frac{\left(\gamma - 3 \right)}{2} \frac{U_3 U_2}{U_1^2} + uB + \left(\gamma -1\right)\frac{U_2^3}{U_1^3} &\frac{\left(3-\gamma \right)}{2} \frac{U_3}{U_1} - \frac{3\left(\gamma -1\right)}{2} \frac{U_2^2}{U_1^2} & \frac{\left(3-\gamma \right)}{2} \frac{U_2}{U_1}& uB'\\
         -\frac{U_4}{U_2}&\frac{-U_4 U_1}{U_2^2}&0&\frac{U_2}{U_1}
        \end{bmatrix}
\end{align}
It can be observed that the flux Jacobian matrix is a function of $\gamma$ and its derivatives given by (\ref{Mixture_gamma_derivatives})
\begin{align}
\label{Mixture_gamma_derivatives}
        B &= \frac{p}{\left(\gamma -1 \right)}\frac{\partial \gamma}{\partial U_1},\\ \nonumber
        B' &= \frac{p}{\left(\gamma -1 \right)}\frac{\partial \gamma}{\partial U_4}.
\end{align}
The flux Jacobian matrix (\ref{Mixture_Flux_Jacobian}) in terms of specific total enthalpy of mixture $H$ is given by (\ref{Mixture_Flux_Jacobian_Enthalpy}).
\begin{align}
\label{Mixture_Flux_Jacobian_Enthalpy}
A(U) = \frac{\partial F\left(U \right)}{\partial U} = \begin{bmatrix}
         0 & 1 & 0& 0 \\ 
         \frac{\left(\gamma-3\right)}{2}u^2 + B & \left(3 - \gamma \right)u &\left(\gamma -1\right) & B'\\ 
          \frac{\left(\gamma -1 \right)}{2}u^3 + Bu -uH & H-\left(\gamma-1  \right)u^2 & \gamma u& B'u\\ 
         -Y_1 u & Y_1&0&u
        \end{bmatrix}
\end{align} 
The eigenvalues of the matrix $A(U)$ \cite{Chargy} are  $(u+a,u,u,u-a$, where $a = \sqrt{\frac{\gamma p}{\rho}})$, and the right eigenvectors are 
\begin{align}
 r_1 =  \begin{bmatrix}
         1 \\ 
         u-a\\ 
         H-ua \\ 
         Y_1
        \end{bmatrix} ,r_2= \begin{bmatrix}
         1 \\ 
         u\\ 
         \frac{u^2}{2}-\frac{B}{\left(\gamma -1\right)}\\
         0
        \end{bmatrix},
r_3 =  \begin{bmatrix}
         0 \\ 
         0\\ 
         -\frac{B'}{\left(\gamma -1\right)} \\
         1
        \end{bmatrix} ,r_4= \begin{bmatrix}
         1 \\ 
         u + a\\ 
         H + ua  \\
         Y_1
        \end{bmatrix}
\end{align}
The above system is hyperbolic as eigenvalues are real and the eigenvectors are linearly independent.
\subsection{Multicomponent fluid simulations with $\gamma$-based model and stiffened gas EOS}
In his fundamental work, Abgrall \cite{RA1} has quoted that any numerical scheme designed for compressible Euler equations extended to multicomponent flows would generate pressure oscillations. It has also been suggested by the author that use of $\frac{1}{(\gamma - 1)}$ as the parameter in the quasi-conservative approach would eliminate the pressure oscillations, this has been demonstrated by Shyue \cite{shyue1,shyue2,shyue3,shyue4}. For $\gamma$ based model the mixture equations in 1D are given by (\ref{gammamodelmixture_governing_equations}).
\begin{align}
\label{gammamodelmixture_governing_equations}
\frac{\partial \rho}{\partial t} + \frac{\left(\partial \rho u\right)}{\partial x} &= 0 \\ \nonumber
\frac{\partial \rho u}{\partial t} + \frac{\partial \left( \rho u^2 + p \right)}{\partial x}&=0\\ \nonumber
\frac{\partial \rho E}{\partial t} + \frac{\partial \left[ \left( \rho E + p \right) u \right]}{\partial x}&=0
\end{align}
Here, $Y^i$ represent the volume fraction of the component gases in a given cell or control volume and $\rho, u, p, p_{\infty}, \gamma$ represent the mixture density, mixture velocity, mixture pressure, mixture stiffened pressure and mixture gamma respectively.
\begin{align}
\label{additional_Equations}
\frac{\left(\partial \frac{\rho}{\left(\gamma - 1\right)} \right)}{\partial t} + \frac{\left(\partial   \frac{\rho u}{\left(\gamma - 1\right)}\right)}{\partial x}=0  \\
\frac{\left(\partial \frac{\rho \gamma p_{\infty}}{\left(\gamma - 1\right)} \right)}{\partial t} + \frac{\left(\partial   \frac{\rho \gamma p_{\infty}u}{\left(\gamma - 1\right)}\right)}{\partial x}=0  \\
\frac{p + \gamma p_{\infty}}{\gamma - 1} = \rho e
\label{stiffenedgaseos}
\end{align}
In the above equations \eqref{stiffenedgaseos} refer to stiffened gas EOS. This EOS would revert to perfect gas EOS in the limit $p_{\infty}\rightarrow0$. These equations \eqref{additional_Equations} along with \eqref{gammamodelmixture_governing_equations} are referred to as $\gamma$-based model by Abgrall \cite{Abgrall_1}. 
\subsection{Numerical Methodology and Novel Algorithms} 
The governing equations are Euler equations in conservative form as given in (\ref{EulerCompact})
\begin{subequations}
\begin{align}
\frac{d \overline{U}}{dt} = - R, 
R &= \frac{1}{\Omega} \left[\sum_{i=1}^{N}{F_c \cdot \hat{n} ~dS}\right], \label{EulerCompact}\\
 \overline{U} &= \frac{1}{\Omega}\int_{\Omega} U d\Omega~.\label{averageU}
\end{align}
\end{subequations}where $U$ is conserved variable vector, $F_c$ is convective flux vector on an interface, $R$ representing net flux from a given control volume, $\Omega$ volume of control volume and $N$ representing number of control surfaces for a given control volume.
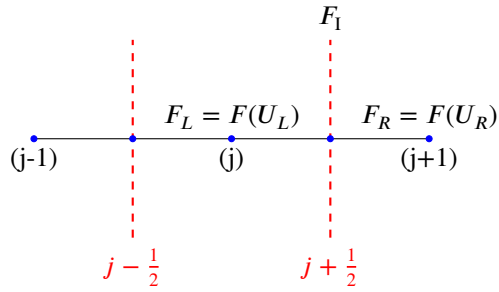
\begin{figure}[ht!]
\begin{center}
\begin{tikzpicture}[scale = 1.3]
\draw (0,1)node (xaxis) [below] {(j-1)}  -- (2,1)node (xaxis) [below] {(j)} -- (4,1)  node (xaxis) [below] {(j+1)};
\draw (2,1) node (xaxis)[above] {$F_L = F(U_L)$};
\draw (4,1) node (xaxis)[above] {$F_R = F(U_R)$};
\draw [red,thick,dashed](1,0) node (yaxis) [below] {$j-\frac{1}{2}$}-- (1,2);
\draw [red,thick,dashed](3,0) node (yaxis) [below] {$j+\frac{1}{2}$}-- (3,2);
\draw(3,2) node (yaxis) [above] {$F_\mathrm{I}$};
\draw [blue,fill] ( 1,1) circle [radius=0.03] ;
\draw [blue,fill] ( 3,1) circle [radius=0.03] ;
\draw [blue,fill] ( 0,1) circle [radius=0.03] ;
\draw [blue,fill] ( 2,1) circle [radius=0.03] ;
\draw [blue,fill] ( 4,1) circle [radius=0.03] ;
\end{tikzpicture}    
\caption{Typical finite volume in 1D}
\label{Interfaceflux1}
\end{center}
\end{figure}
The convective flux on any interface of a control volume as shown in figure(\ref{Interfaceflux1}) for any stable scheme can be written as the sum of an average flux across the interface and a numerical dissipative flux as given in \eqref{ConvectiveDiffusiveFlux}.
\begin{subequations}
\begin{align}
 F_\mathrm{I} = \frac{1}{2}\left[ F_L + F_R\right] - d_\mathrm{I} \label{ConvectiveDiffusiveFlux} \\
d_\mathrm{I}  = \frac{\mid \alpha_{\mathrm{I}} \mid}{2}\left(U_R - U_L\right) \label{DiffusionFlux}
\end{align}
\end{subequations}
where $\alpha_{\mathrm{I}}$ coefficient of numerical diffusion and $d_\mathrm{I}$ represents the numerical dissipative flux. In the present work the coefficient of numerical diffusion is determined by MOVERS \cite{Jaisankar_SVRRao}, RICCA and MOVERS+ \cite{rkolluru} as briefed in sections(\ref{sec:MOVERS},\ref{sec:RICCA} $\&$ \ref{sec:movers+}) are utilised to simulate the mixture equations based on mass fraction model with perfect gas EOS and $\gamma$ based model with stiffened gase EOS. These algorithms are independent of eigenstructure of the underlying hyperbolic system and can be extended to any arbitrary EOS.

\subsubsection{MOVERS} \label{sec:MOVERS}
The central scheme of interest is due to \cite{Jaisankar_SVRRao} who introduced a new central scheme named MOVERS (\emph{Method of Optimal Viscosity for Enhanced Resolution of Shocks}) which can capture grid aligned shocks and contact-discontinuities accurately. This coefficient of numerical diffusion rewritten in terms of RH conditions is given by \eqref{alphai}
\begin{align}
\label{alphai}
 \lvert \alpha_{I} \rvert_i = \left| s_i \right| =\left| \frac{\Delta F_i}{\Delta U_i} \right|, \quad i =1,2,3, \quad \Delta(\cdot) = (\cdot)_R - (\cdot)_L
\end{align}
In order to introduce boundedness and further stabilize the numerical scheme, $\alpha_I$ is to be restricted to a physically feasible range of eigenvalues of the flux Jacobian matrix. This process known as wave speed correction \eqref{wavespeedcorrection} is incorporated such that the coefficient of numerical diffusion lies within the eigenspecturm of the flux Jacobian \textit{i.e.}, $\alpha_I \in 
\left[\lambda_{max}, \lambda_{min}\right]$.
\begin{align}
\label{wavespeedcorrection}
 \lvert \alpha_{\mathrm{I}} \rvert = \begin{cases}
                             \lambda_{\textit{max}},~~ if~~ \lvert \alpha_{\mathrm{I}} \rvert > \lambda_{\textit{max}} \\
                             \lambda_{\textit{min}},~~ if~~ \lvert \alpha_{\mathrm{I}} \rvert < \lambda_{\textit{min}}\\
                             \lvert \alpha_{\mathrm{I}} \rvert , \\
                            \end{cases}
\end{align}
This method is independent of eigenstructure of the underlying hyperbolic systems, is simple and can capture grid-aligned stationary discontinuities exactly. Authors \cite{Jaisankar_SVRRao} introduced two variations of MOVERS: $(i)$ an $n$-wave based coefficient of numerical diffusion, corresponding to $n$ number of conservation laws (MOVERS-n) and $(ii)$ a scalar diffusion, corresponding to the energy equation, referred to as MOVERS-1. 
The robustness of the basic scheme has been improvised through its variants by Maruthi N.H. \cite{Maruthi_Thesis} and extended them to hyperbolic systems for magnetohydrodynamics and shallow water flows. In this work this algorithm is chosen as the foundation to devise two new efficient algorithms for hyperbolic systems. The simplicity and accuracy of this algorithm make this scheme a well-suited base-line solver for further research, apart from its independency of the eigenstrucure.

\subsubsection{Riemann Invariant based Contact-discontinuity Capturing Algorithm (RICCA)}\label{sec:RICCA}
 
The numerical diffusion evaluated using Riemann Invariant based Contact-discontinuity Capturing Algorithm (RICCA) is given by
\begin{equation}\label{accudisks_eval_alpha_euler_2d}
{\alpha}_ {\emph{$I$}} = \begin{cases}
                            \qquad \qquad \frac{|V_{nL}|+|V_{nR}|}{2}, \qquad \qquad \qquad \quad \text{if } |\Delta\mathbf{F}|<\delta \ \text{and} \ |\Delta\mathbf{U}| <\delta\\
                            \quad \qquad \qquad \qquad \qquad \qquad \qquad \qquad \qquad  \\
                            max(|V_{nL}|, |V_{nR}|) + sign(|\Delta p_{\raisebox{-2pt} {\scriptsize {$\mathrm{I}$}}}|) a_{\mathrm{I}} , \qquad \text{otherwise} \\
                            \qquad \qquad \qquad \qquad \quad
                           \end{cases}
\end{equation}
where $a_{\mathrm{I}} = \sqrt{\frac{\gamma p_{\raisebox{-2pt} {\scriptsize \emph{$I$}}}}{\rho_{\raisebox{-2pt} {\scriptsize \emph{$I$}}}}}$ is the speed of sound evaluated with the values at the interface given by
\begin{align}
\rho_{\raisebox{-2pt} {\scriptsize \emph{$I$}}} = \frac{\rho_L+\rho_R}{2}, 
p_{\raisebox{-2pt} {\scriptsize \emph{$I$}}} = \frac{p_L+p_R}{2},
\Delta{p}_{\raisebox{-2pt} {\scriptsize \emph{$I$}}} = (p_R-p_L).
\end{align}

\subsubsection{MOVERS without wave speed correction - \emph{MOVERS+}} \label{sec:movers+}
The coefficient of numerical diffusion for MOVERS+ is given by
\begin{align}
\label{MOVERS_NWSC}
 \lvert d_{\mathrm{I}} \rvert_j =  \Phi Sign(\Delta U_j)\lvert \Delta F_j \rvert + \left(\frac{|V_{nL}|+|V_{nR}|}{2}\right) \Delta U_j, \quad j =1,2,3 
\end{align}
These two new algorithms RICCA and MOVERS+ 
\begin{itemize}
 \item can capture steady contact-discontinuities exactly,
 \item has sufficient numerical diffusion near shocks so as to avoid shock instabilities, and
 \item does not need entropy fix for at sonic points.
\end{itemize}
A similar strategy was introduced by N.Venkata Raghavendra in \cite{Venkat_thesis,Venkat_Arxiv} to design an accurate contact-discontinuity capturing discrete velocity Boltzmann scheme for inviscid compressible flows.

\subsection{Modifications of upwind methods for multicomponent flows}
As mentioned before, the application of upwind methods to multicomponent flows is non-trivial because these methods are strongly dependent on the eigenstructure. Larrouturou and Fezoui \cite{Larrouturou1} have reviewed these modifications needed for upwind methods, which are briefly presented here.
\subsubsection{Extension of Steger-Warming FVS method to multicomponent gases}
The Flux Vector Splitting (FVS)  method of Steger-Warming method as given in \cite{TORO_1}, leads to the following split flux vectors for Euler equations. 
\begin{align}
\label{stegerwarmingfluxes}
F^{\pm} = \frac{\rho}{2\gamma} \begin{bmatrix}
         \left(u-a\right)\lambda_1^{\pm} + 2\left(\gamma -1 \right)\lambda_2^{\pm} + \lambda_3^{\pm} \\ 
         \lambda_1^{\pm} + 2\left(\gamma -1 \right)u\lambda_2^{\pm} + \left(u+a\right)\lambda_3^{\pm} \\ 
         \left(H-ua\right)\lambda_1^{\pm} + 2\left(\gamma -1 \right)u^2\lambda_2^{\pm} + \left(H + ua\right)\lambda_3^{\pm} \\ 
          \end{bmatrix} \\ \nonumber
\end{align}
 For the mixture equations (\ref{mixture_governing_equations}), an additional fourth component (for the extra mass fraction term) for the split flux vectors is given by 
 \begin{align}
  F[4] = \rho u Y_1 = F^{\pm}[1]Y_1.
  \end{align}
\subsubsection{Extension of van Leer FVS method to multicomponent gases}
The details of the flux vector splitting developed by van Leer are given in \cite{TORO_1} for Euler equations. 
Extension of van Leer flux splitting to multicomponent mixture equations as a function of Mach number is given 
by 
\begin{align}
 F = F(\rho,a,M,Y) =  \begin{bmatrix}
         \rho a M \\ 
         \rho a^2 \left( M^2 + 1 \right)\\ 
         \rho a^3 M \left( \frac{M^2}{2} + \frac{1}{\left(\gamma -1\right)} \right) \\ 
         \rho a M Y_1
        \end{bmatrix}
\end{align}
The split fluxes given in \cite{Larrouturou1} are 
\begin{align}
 F^{\pm} = \frac{1}{4} \rho a \left(1{\pm} M \right)^2 \begin{bmatrix}
         1 \\ 
         \frac{2a}{\gamma}\left(\frac{\left(\gamma -1\right)}{2}M {\pm} 1\right)\\ 
          \frac{2a^2}{\left(\gamma^2-1\right)}\left(\frac{\left(\gamma -1\right)}{2}M {\pm} 1\right)^2\\ 
          Y_1
        \end{bmatrix}
\end{align}
\subsubsection{Extension of Roe's FDS method}
Roe's Flux Difference Splitting (FDS) method, which is an approximate Riemann solver, cannot be directly extended to multicomponent flows in a trivial way. 
 In order to evaluate Roe's numerical flux, the following information is necessary
 \begin{itemize}
  \item wave strengths $\tilde{\alpha_i}$,
  \item eigenvalues of the flux Jacobian matrix $\tilde{\lambda_i}$,
  \item right eigenvectors of the flux Jacobian matrix $\tilde{R^{\left(i\right)}}$.
 \end{itemize}
 The following basic conditions (also called as $U$ property) are to be satisfied by Roe scheme  
\begin{enumerate}
 \item consistency, $A\left(U_L,U_R\right) = A(U)$ if $U_L=U_R=U$,
 \item hyperbolicity \ie flux Jacobian matrix should have real eigenvalues,
 \item conservation across discontinuities, $F\left(U_R\right) - F\left(U_L\right) = A\left(U_R - U_L\right)$.
\end{enumerate}
 The two component scheme for the interface flux given in \cite{Larrouturou1} is
 \begin{equation}
  F\left(U_L,U_R \right) = \frac{1}{2}\left( F(U_L) + F(U_R) \right) + \frac{1}{2}|\tilde{A}|\left(U_L-U_R\right)
 \end{equation}
where $U$ represents the average state between the left and right states. The averaged state is defined as 
\begin{align}
\begin{rcases}
 U &= \left[\tilde{\rho}, \tilde{\rho} \tilde{u},\tilde{\rho}\tilde{E},\tilde{\rho}\tilde{Y} \right]^T \\
 \tilde{\rho} &= \frac{\rho_L\sqrt{\rho_L} + \rho_R\sqrt{\rho_R}}{\sqrt{\rho_L} + \sqrt{\rho_R}}\\
 \tilde{u} &= \frac{u_L\sqrt{\rho_L} + u_R\sqrt{\rho_R}}{\sqrt{\rho_L} + \sqrt{\rho_R}}\\
 \tilde{H} &= \frac{H_L\sqrt{\rho_L} + H_R\sqrt{\rho_R}}{\sqrt{\rho_L} + \sqrt{\rho_R}}\\
 \tilde{Y} &= \frac{Y_L\sqrt{\rho_L} + Y_R\sqrt{\rho_R}}{\sqrt{\rho_L} + \sqrt{\rho_R}}\\
 \tilde{A}(U) &= \begin{bmatrix}
         0 & 1 & 0& 0 \\
         \frac{\left(\tilde{\gamma}-3\right)}{2}\tilde{u}^2 + \tilde{B} & \left(3 - \tilde{\gamma} \right)\tilde{u} &\left(\tilde{\gamma} -1\right) & B'\\ 
          \frac{\left(\tilde{\gamma} -1 \right)}{2}\tilde{u}^3 + \tilde{B}\tilde{u} -\tilde{u}\tilde{H} & \tilde{H}-\left(\tilde{\gamma}-1  \right)\tilde{u}^2 & 
\tilde{\gamma}\tilde{u}& \tilde{B}'\tilde{u}\\ 
         -\tilde{Y_1} \tilde{u} & \tilde{Y_1}&0&\tilde{u}
        \end{bmatrix}
        \label{roematrix}
        \end{rcases}
\end{align}
The matrix $\tilde{A}$ as given in (\ref{roematrix}) is diagonalisable, and its eigenvalues, as given in \cite{Chargy}, are $\left(\tilde{u} - \tilde{a},\tilde{u} ,\tilde{u},\tilde{u} + 
\tilde{a}\right)$  where $ \left(\tilde{a}^2 = \left(\tilde{\gamma}-1 \right)\left(\tilde{H} - \frac{\tilde{u}^2}{2}\right)\right)$.
\newline
 The following remarks are given in \cite{Chargy,RA4}, which highlights the conditions under which Roe scheme is not applicable.
\begin{remark}
For two component fluid flow, the conservation property is satisfied only if $~\gamma_1~=\gamma_2 ~=\textrm{Constant}$. If $\gamma_1 \neq \gamma_2$, then pressure 
oscillations are observed in the case of steady contact discontinuities.
\end{remark}
\begin{remark}
In order to satisfy the conservation property $\tilde{A}\left(U\right)$ has to be modified to $A\left(\tilde{U}\right)$, for which the following definitions of $\tilde{B}' = 
\frac{c_{v1}c_{v2}\left(\gamma_1 - \gamma_2\right)\tilde{T}}{\tilde{Y_1}c_{v1} + \left(1- Y_1\right)c_{v2}}$ and $\tilde{T} = \frac{T_L\sqrt{\rho_L} + T_R\sqrt{\rho_R}}{\sqrt{\rho_L} + \sqrt{\rho_R}} \neq T\left(\tilde{U}\right)$ are used. The Jacobian matrix $\tilde{A}$ is diagonalisable only if this modification is incorporated.
\end{remark}
\begin{remark}
This expression $\tilde{B}' = 
\frac{c_{v1}c_{v2}\left(\gamma_1 - \gamma_2\right)\tilde{T}}{\tilde{Y_1}c_{v1} + \left(1- Y_1\right)c_{v2}}$ is not easily extendable to a mixture of more than 2 components.
\end{remark}
\begin{remark}
 If the two states $U_L$ and $U_R$ are supersonic and satisfy $u_L \geq a_L$, $u_R\geq a_R$, and when using upwind schemes, the flux $\phi\left(U_L,U_R\right) = 
F\left(U_L\right)$, is satisfied in all upwind schemes, modified for multicomponent cases, except for Roe scheme.
\end{remark}
\begin{remark}
 The construction of Jacobian matrix of average the state $\tilde{A}\left(U_L,U_R\right) = A\left(\tilde{U}\right)$, is impossible for acomplex equation of state for multicomponent cases.
\end{remark}
It can be observed that direct extension of Roe scheme to multicomponent fluids is not an easy task. The new central solvers introduced in section(\ref{sec:MOVERS},\ref{sec:RICCA},\ref{sec:movers+}), RICCA and MOVERS+, along with MOVERS-1 and MOVERS-n, do not require any of the above modifications as they are not dependent on eigenstructure at all. 
\section{Results for mass fraction based model}
In this section the test cases used for the validation of the central solvers, MOVERS-1, MOVERS-n, MOVERS+ and RICCA, for multicomponent flows are discussed. Initial conditions for the 1-D shock tube are given in table (\ref{Test_Cases_1}).

\begin{table}[h!]
\centering
\begin{tabular}{ |c|c|c|c|c|c|c|c|c|c|}
\hline
Test case & $\rho_L$    &    $H_L$     &     $u_L$     &     $\rho_R    $     &     $H_R$         &     $u_R$     & $\gamma_L$& $\gamma_R$&     Time \\
\hline
1&        1.0    &    1.0    &     0.0         &    0.125    &    1.0            &     0.0         & 1.4 &1.6&steady case    \\
\hline 
\end{tabular}
\caption{Steady contact-discontinuity test case with perfect gas EOS}
\label{Steady_CD_Testcase}
\end{table}
\begin{table}[h!]
\centering
\begin{tabular}{ |c|c|c|c|c|c|c|c|c|c|}
\hline
Test case & $\rho_L$    &    $p_L$     &     $u_L$     &     $\rho_R    $     &     $p_R$         &     $u_R$     & $\gamma_L$& $\gamma_R$&     Time \\
\hline
1&        1.0    &    1.0    &     -1.0         &    1.0    &    5.0            &     1.0         & 1.4 &1.4&0.21    \\
\hline
2&        1.0    &    1000.0    &     0.0         &    0.125    &    1.0            &     0.0         & 1.6 &1.4&0.21    \\
\hline 
\end{tabular}

\caption{Initial conditions for shock tube problem and mass fraction positivity test case:  data obtained from \cite{karni1}, \cite{Chargy} with Perfect gas EOS}
\label{Test_Cases_1}
\end{table}

 A fully coupled approach is used for the flux evaluation approach as explained in \cite{Chargy} and it is stated that the Steger-Warming scheme and van Leer Scheme will preserve the maximum of $0 \leq Y \leq 1$ as they are using a fully coupled approach. 

\subsection{Steady contact-discontinuity}
This test case refers to a contact discontinuity wherein there is a jump in density and $\gamma$ as given in table (\ref{Steady_CD_Testcase}). The shock tube is filled with two different perfect gases denoted by variable $\gamma$. The initial discontinuity is present at $x=0.5$ with $x\in [0,1]$. A total of 100 control volumes are used in simulations.  Various numerical schemes like Steger Warming, VanLeer, MOVERS-n,MOVERS-1,RICCA and MOVERS+ are compared for this test case. The ability of the numerical schemes to resolve the steady contact-discontinuity is analysed here.
\begin{figure}[htb]
\begin{center}
 \includegraphics[scale=0.7]{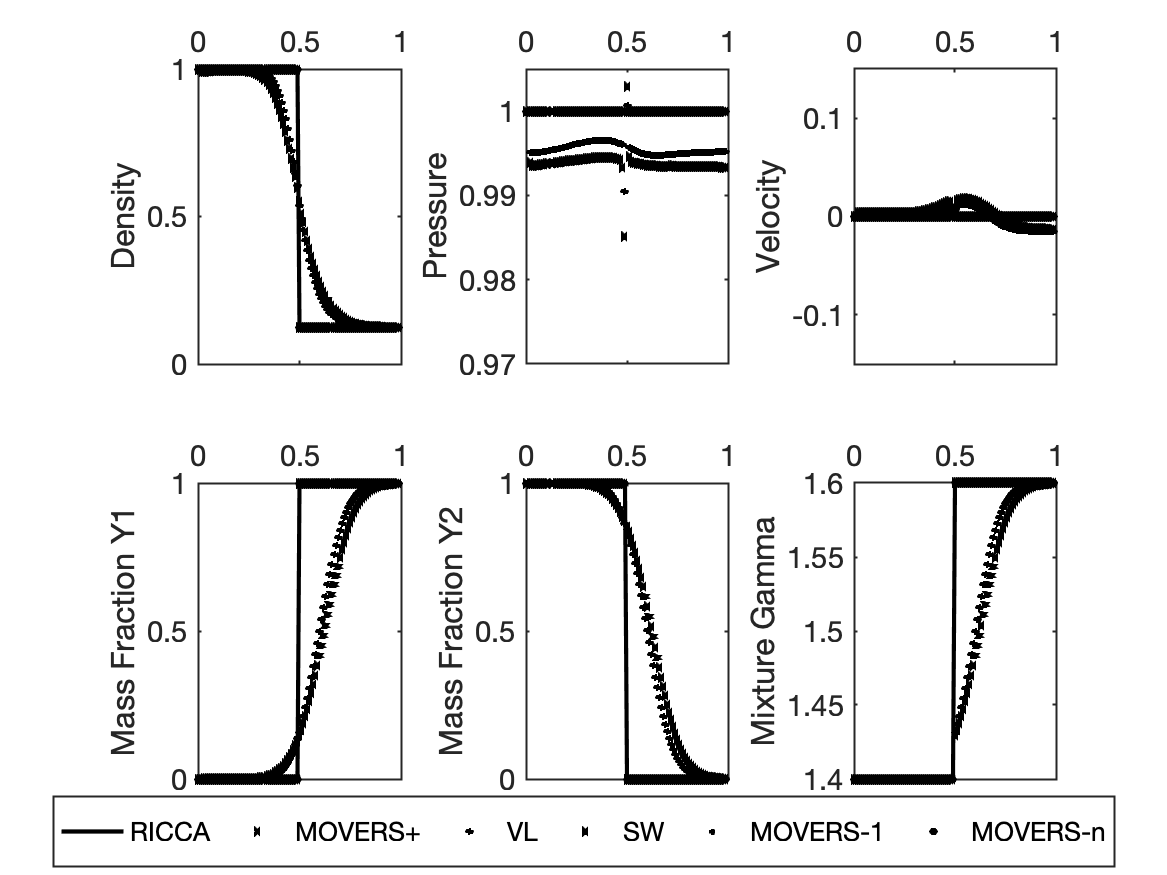} 
\caption{Steady contact-discontinuity Comparison of RICCA,MOVERS+,MOVERS-n, MOVERS-1,Vanleer and Steger Warming with unequal $\gamma$}
\label{Comparision_Steady_Contact}
\end{center}
\end{figure}

Figures (\ref{Comparision_Steady_Contact}) refer to comparison of solution obtained for the steady contact-discontinuity case as given in table (\ref{Steady_CD_Testcase}). It can be observed that RICCA, MOVERS+, MOVERS-1 and MOVERS-n resolve the steady contact exactly and even the mass fraction is resolved exactly, whereas for Steger-Warming and van Leer methods, the contact discontinuity and the mass fractions are diffused. It can also be observed that in Steger-Warming scheme and van Leer scheme oscillations are present in pressure and velocity but the positivity of the mass fraction is preserved.

\subsection{Pressure oscillations test case}
This test case (\ref{Test_Cases_1}) is used by \cite{karni2}, to test the positivity of the mass fraction by the regular Godunov type conservative finite volume methods. In \cite{karni2}, the authors claim that many of the numerical methods which are formulated in the conservative finite volume method would fail to preserve the positivity of the mass fraction. Further, the pressure and velocity have oscillations for the regular finite volume methods.
\begin{figure}[htb]
\begin{center}
 \includegraphics[scale=0.7]{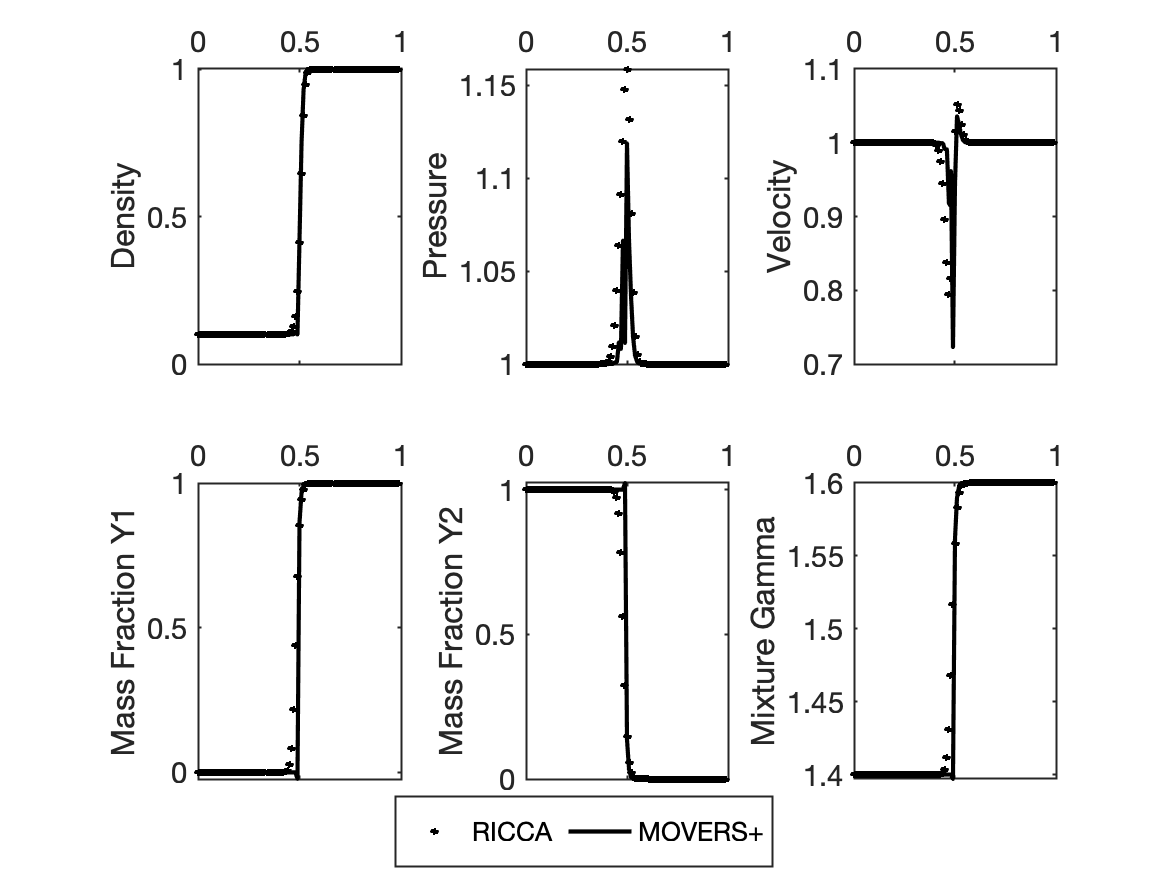}
\caption{Isolated material front problem using RICCA and MOVERS+}
\label{MOVERSP_Isolated_Material_Front}
\end{center}
\end{figure}
Figures ( \ref{MOVERSP_Isolated_Material_Front}) refer to isolated material front test case as described in table (\ref{Test_Cases_1}). As commented by the authors in \cite{karni3}, for this test case pressure oscillations are present for all first-order numerical schemes which are designed based on the conservative formulation with fully coupled approach and the positivity of the mass fractions for such scheme is doubtful. Though the new algorithms are based on a fully coupled approach in conservative formulation with controlled numerical diffusion, they produce mild oscillations in pressure but they preserve the positivity of the mass fraction unlike the other conservative schemes. 
\subsection{Sod shock tube test case}
This is a standard test case whose initial conditions are given in table (\ref{Test_Cases_1}) for variable gamma values. The data values of the test case are taken from \cite{karni1} also referred in \cite{Chargy}. The second test case is a stiffer shock tube problem with variation of pressure, mentioned as in test case 3 in table (\ref{Test_Cases_1}). For these shock tube problems the initial discontinuity is located at $x=0.5$, with mass fraction to the left of the discontinuity $Y_L =1, Y_R=0$ and to the right of the discontinuity $Y_L=0,Y_R =1$. In all of the computations 100  equally spaced control volumes on the interval $[0,1]$ are considered with CFL of $0.45$ till the prescribed time is reached. Numerical results are presented for Steger-Warming scheme, van Leer scheme, MOVERS-1, MOVERS-n, MOVERS+, and RICCA.

\begin{figure}[htb]
\begin{center}
\includegraphics[scale = 0.7]{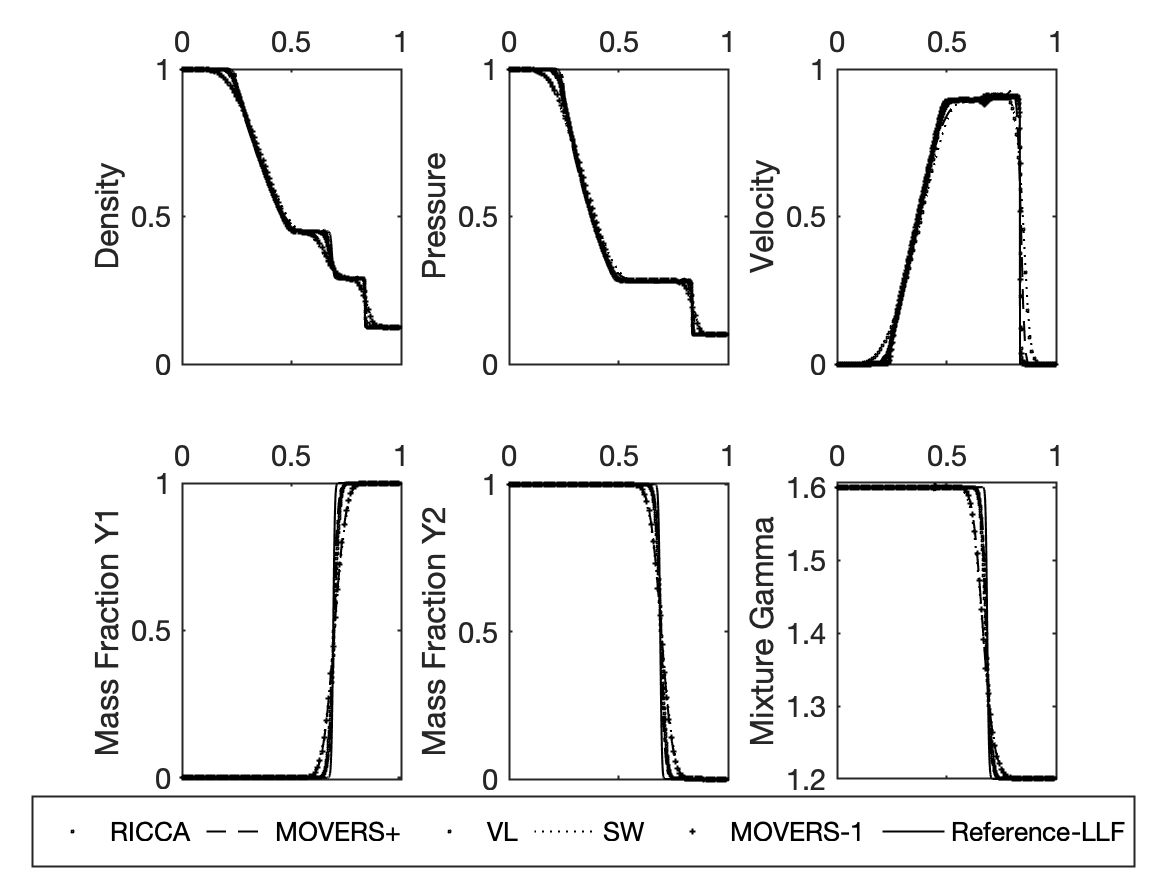} 
\caption{Sod shock tube problem using RICCA with unequal $\gamma$, $\gamma_L = 1.6$, $\gamma_R = 1.2$}
\label{All_SOD_Unequalgamma}
\end{center}
\end{figure}

\begin{figure}[htb]
\begin{center}
            \includegraphics[scale=0.7]{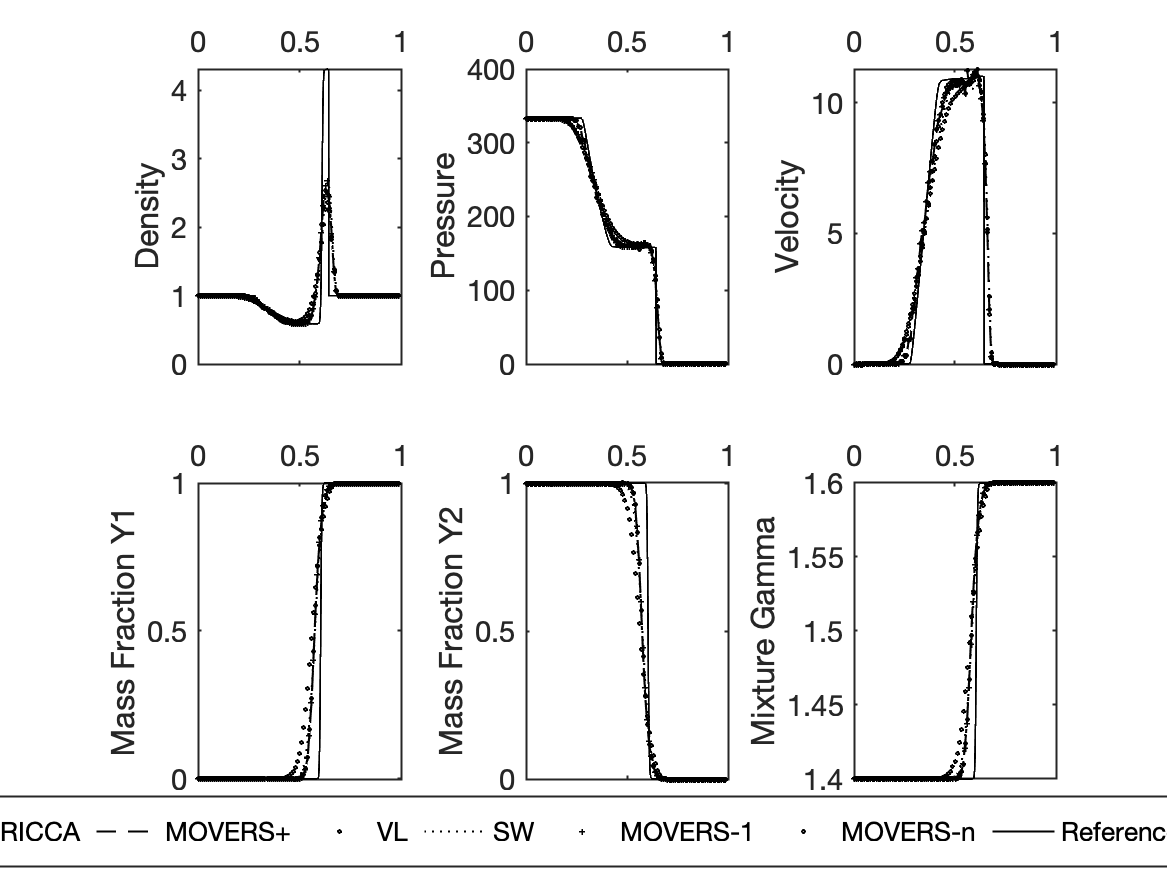}
\caption{Stiff shock tube problem using RICCA,MOVERS+, MOVERS-n,MOVERS-1 with equal gamma $\gamma = 1.4$}
\label{All_StiffShockTube_equalgamma}
\end{center}
\end{figure}

Figures (\ref{All_SOD_Unequalgamma}) refer to the standard Sod shock tube problem whose initial conditions are defined as test case 2 in table (\ref{Test_Cases_1}). This shock tube problem has two different fluids with different $\gamma$ values initially separated by a diaphragm placed at $x=0.5$. Numerical simulations are carried out with 100 control volumes for all schemes and the reference solution is generated using 10000 control volumes using Rusanov Method. 
Figures ( \ref{All_StiffShockTube_equalgamma}) refer to the 
stiff shock tube whose initial conditions are defined as in test case 3 in table (\ref{Test_Cases_1}). As can be seen all the numerical schemes preserve the mass fraction positivity and no pressure oscillations are present in the pressure.

\subsection{1D results for $\gamma$-based model}
Numerical simulations have been carried out for the test cases described in \cite{shyue1} whose initial conditions are given in the table (\ref{table:MultiComponent_Stiffened_gas_EOS_Test_Cases}).
\begin{table}[h!]
\centering
\begin{tabular}{ |c|c|c|c|c|c|c|c|c|c|c|c|}
\hline
Sno & $\rho_L$    &    $p_L$     &     $u_L$     &     $\rho_R    $     &     $p_R$         &     $u_R$     & $\gamma_L$& $\gamma_R$& ${p_{\infty}}_L$ &${p_{\infty}}_R$&     Time \\
\hline
1&        1.0    &    1.0    &     1.0         &    0.125    &    1.0            &     1.0         & 1.4 &1.2&0&0&0.12    \\
\hline 
2&        1.0    &    1.0    &     1.0         &    0.125    &    1.0            &     1.0         & 1.4 &4.0&0&1&0.12    \\
\hline 
3&        1.241    &    2.753    &     0.0         &    1.0    &    $3.059 \times 10^{-4}$            &     0.0         & 1.4 &5.5&0&1.505&0.1    \\
\hline 
4&        1.0    &    1.0    &     0.0         &    5.0    &    1.0            &     0.0         & 1.4 &4.0&0&1&0.2    \\
\hline 
&            &        &             &    7.093    &    10.0            &     -0.7288         &  &4.0&0&1&0.2    \\
\hline 
\end{tabular}
\caption{1D Shock tube test cases referred from \cite{shyue1} with stiffened gas EOS}
\label{table:MultiComponent_Stiffened_gas_EOS_Test_Cases}
\end{table}

\begin{figure}[htb]
\begin{center}
            \includegraphics[scale=0.75]{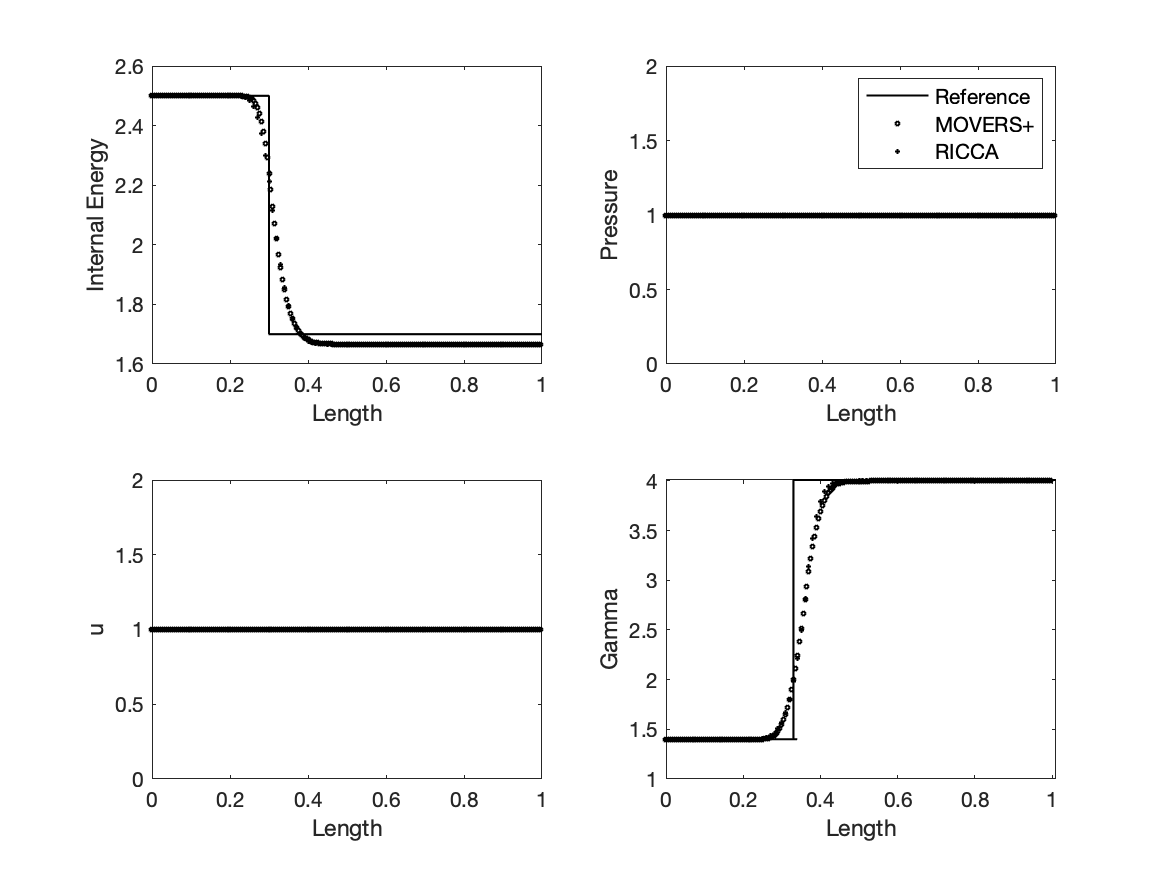}
            \caption{Interface only problem with stiff gas simulated using RICCA using eq(10) of \cite{shyue1}}
	\label{InterfaceOnlyStiffgaswitheq10RICCA}
\end{center}
\end{figure}

Test case 1 is an interface only problem and consists of a single contact discontinuity. This test case consists of two sets of data 
\begin{enumerate}
\item a polytropic gas with two constants states as case 1,
\item has same states except for the changes in $\gamma$ and $P_{\infty}$ as in case 2.
\end{enumerate}
Initial position of the diaphragm is located at $x = 0.2$ and the length of the shock tube, $L=1$. Results are shown for 1O RICCA and MOVERS+. 100 control volumes are considered for computation for both cases and the results are shown at the prescribed time of 0.12. 
\begin{figure}[htb]
\begin{center}
            \includegraphics[scale=0.75]{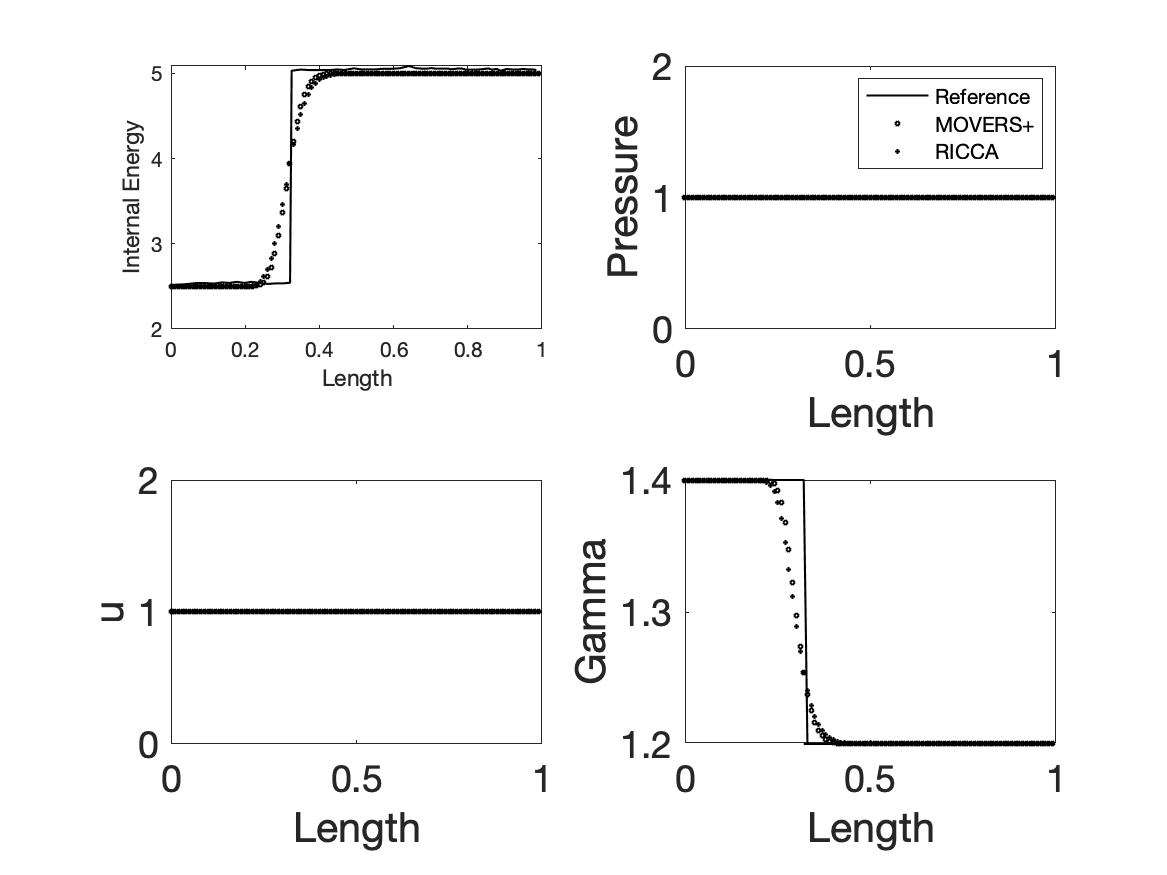}
            \caption{Interface only problem with $P_{\infty} = 0$ simulated using MOVERS+ and RICCA using eq(10) in \cite{shyue1}}
\label{InterfaceOnlyNoStiffgaswitheq10RICCA}            
\end{center}
\end{figure}
  Figures (\ref{InterfaceOnlyNoStiffgaswitheq10RICCA},\ref{InterfaceOnlyNoStiffgaswitheq10NWSC}) refer to interface only problem with $p_{\infty} =0$, which corresponds to perfect gas EOS using RICCA and MOVERS+. It can be observed from the figures that the pressure oscillations are not present when gamma based model (unlike in Figure (\ref{gammamodelmixture_governing_equations})) is used and $\frac{1}{\gamma -1}$ is used as the conservative variable as suggested by \cite{RA1} and reconfirmed by \cite{shyue1}. Shyue has also suggested that using eq(10) in \cite{shyue1} cannot be generalised to all shock interaction problems. 
Figures (\ref{InterfaceOnlyStiffgaswitheq10RICCA},\ref{InterfaceOnlyStiffgaswitheq10NWSC}) refer to interface only problem with $P_{\infty} \neq 0$, which corresponds to stiffened gas EOS using RICCA and MOVERS+ schemes. It can be observed from the figure that the pressure oscillations are not present when gamma based model. Even for a large jump in $\gamma$ both the numerical schemes do not generate any pressure oscillations.

\begin{figure}[htb]
\begin{center}
            \includegraphics[scale=0.75,keepaspectratio]{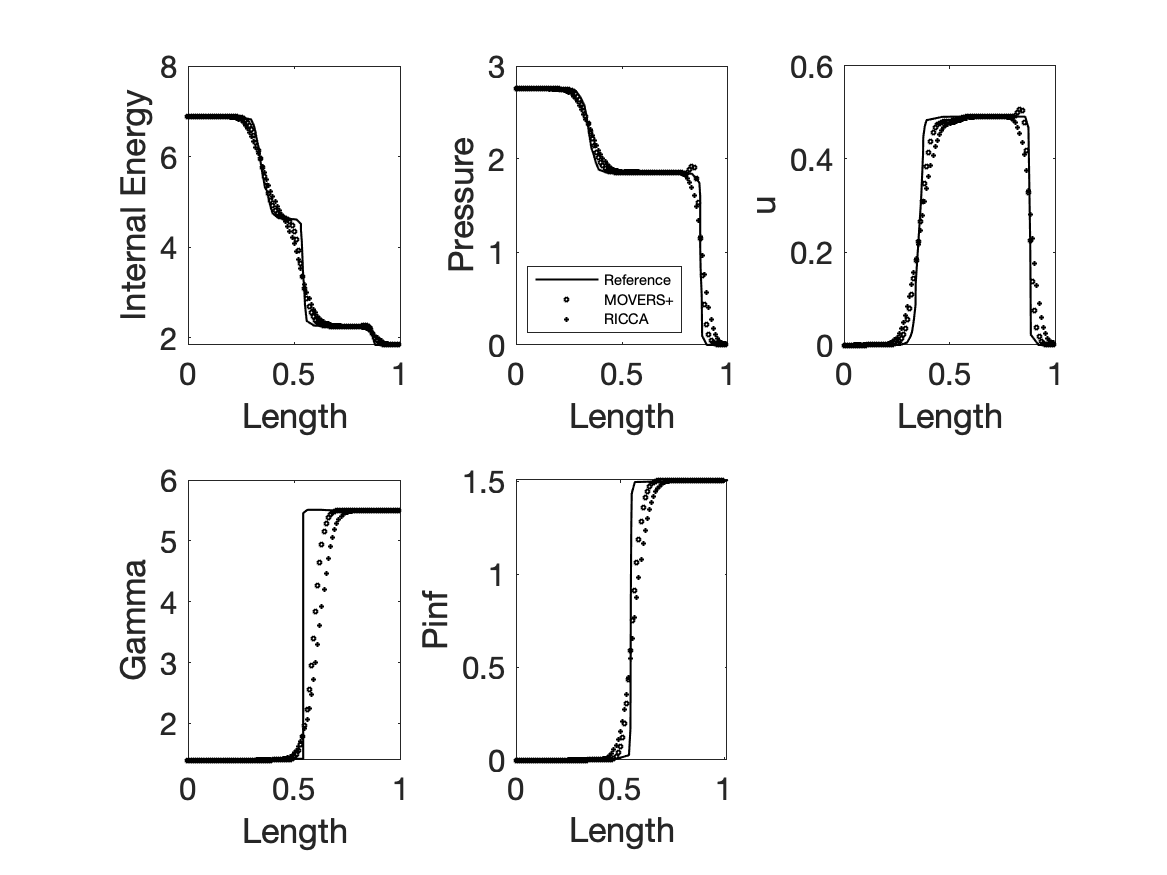}
            \caption{Liquid gas Riemann problem using MOVERS+ and RICCA}
\label{Liquid_Gas_RP_NWSC}
\end{center}
\end{figure}

Test case 2 is a two fluid gas-liquid Riemann problem with initial conditions as given in case 3 in the table. The diaphragm position is located with gas occupying the domain till $x\leq0.5$ and then the liquid and the time for computation is $t = 0.1$. In these test cases the conservative formulation using  $\frac{\rho}{\gamma - 1}$ is used. Data for reference solution is taken from \cite{shyue1}. Simulations are shown for MOVERS+ and RICCA with100 control volumes. It can be observed form the figures (\ref{Liquid_Gas_RP_RICCA},\ref{Liquid_Gas_RP_NWSC}) and figure (\ref{Liquid_Gas_RP_Pinf}) that there are no pressure oscillations present in the pressure and the velocity. The internal energy and the $p_{\infty}$ are accurately predicted.

\begin{figure}[htb]
\begin{center}
            \includegraphics[scale=0.75]{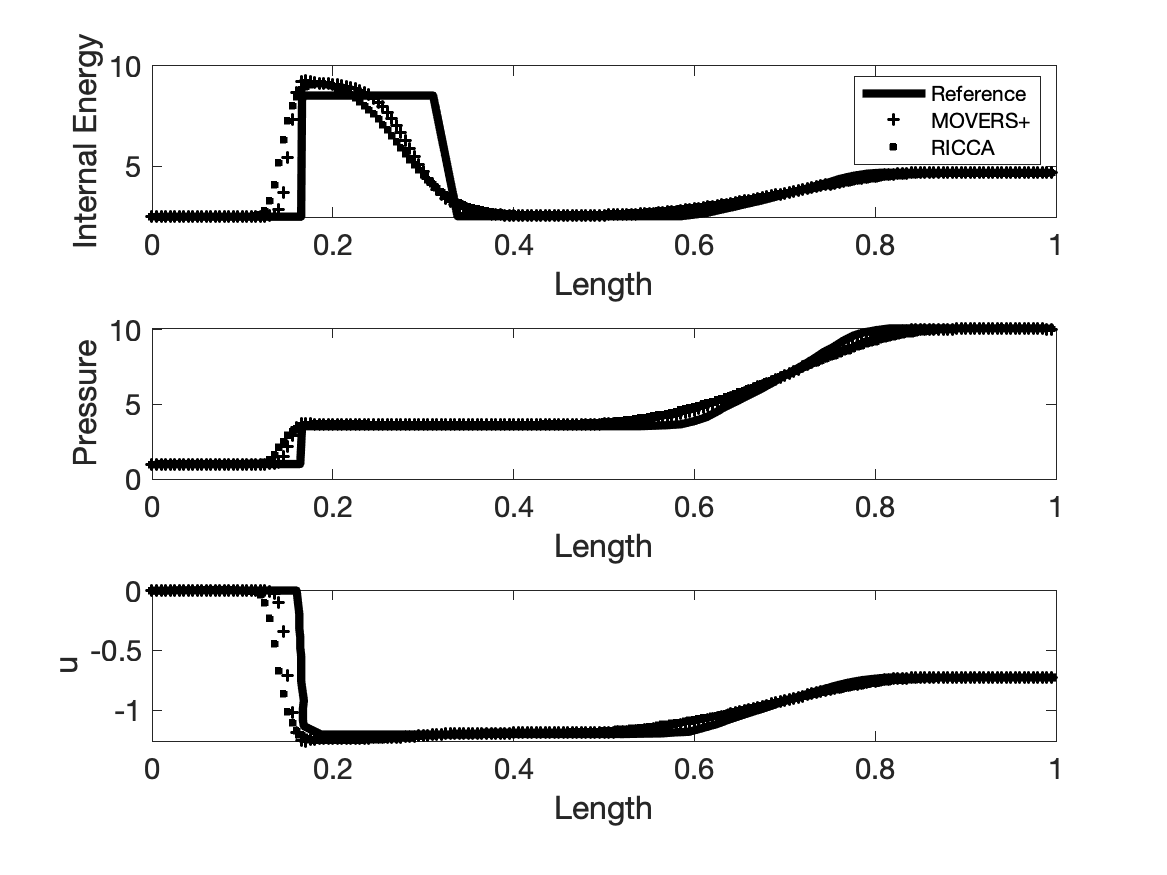}
            \caption{Sock contact interaction problem using RICCA and MOVERS+}
\label{SCI_NWSC}
\end{center}
\end{figure}

The third test case considered here is a shock contact-discontinuity interaction problem with the data given as in case 4. Here the two liquids are separated by the interface at $x=0.5$ and a shock wave located at $x= 0.6$ with the pre- and post-shock conditions as given in the table \ref{table:MultiComponent_Stiffened_gas_EOS_Test_Cases} and the computations are carried out for a time $t=0.2$. For this case 200 control volumes are considered in the computation, the results from computation are shown in figure (\ref{SCI_NWSC}). It can be seen that the phenomenon is captured accurately by both RICCA and MOVERS+. 

\subsection{Extension to two dimensions}
In this section, the multicomponent model described in the previous section is extended to 2D. The governing equations for $\gamma$-based model in 2D are given by (\ref{2DMulticompoent_geqns})
\begin{align}
\label{2DMulticompoent_geqns}
\frac{\partial \rho}{\partial t} + \frac{\left(\partial \rho u\right)}{\partial x} + \frac{\left(\partial \rho v\right)}{\partial y} &= 0 \\
\frac{\partial \rho u}{\partial t} + \frac{\partial \left( \rho u^2 + p \right)}{\partial x}+ \frac{\left(\partial \rho u v\right)}{\partial y}&=0\\
\frac{\partial \rho v}{\partial t} + \frac{\left(\partial \rho u v\right)}{\partial x} + \frac{\partial \left( \rho v^2 + p \right)}{\partial y}&=0\\
\frac{\partial \rho E_t}{\partial t} + \frac{\partial \left[ \left( \rho E_t + p \right) u \right]}{\partial x}+ \frac{\partial \left[ \left( \rho E_t + p \right) v \right]}{\partial y}&=0
\end{align}
with the equations $\gamma$ given by 
\begin{align}
\frac{\left(\partial \frac{\rho}{\left(\gamma - 1\right)} \right)}{\partial t} + \frac{\left(\partial   \frac{\rho u}{\left(\gamma - 1\right)}\right)}{\partial x}+ \frac{\left(\partial   \frac{\rho v}{\left(\gamma - 1\right)}\right)}{\partial x}&=0  \\
\frac{\left(\partial \frac{\rho \gamma p_{\infty}}{\left(\gamma - 1\right)} \right)}{\partial t} + \frac{\left(\partial   \frac{\rho \gamma P_{\infty}u}{\left(\gamma - 1\right)}\right)}{\partial x}+ \frac{\left(\partial   \frac{\rho \gamma p_{\infty}v}{\left(\gamma - 1\right)}\right)}{\partial y}&=0  
\end{align}
and stiffened gas EOS as given by (\ref{2Dstiffenedgaseos_additional_Equations}).
\begin{align}
\frac{p + \gamma p_{\infty}}{\gamma - 1} = \rho e
\label{2Dstiffenedgaseos_additional_Equations}
\end{align}

For the 2D $\gamma$-based model given above, simulations are carried out using RICCA and MOVERS+.
\subsubsection{Moving Interface Problem}
The first test case considered is a moving interface problem which consists of a bubble with radius $r_0 = 0.16$ evolving in a constant velocity field $\left( u, v \right) = (1,1)$. The initial data considered here is similar to the 1D test case as described in table (\ref{table:MultiComponent_Stiffened_gas_EOS_Test_Cases}). The pressure is uniform with value $p=1$ while the $\rho,\gamma ,p_{\infty}$ jump across the interface. Initially the bubble is placed at $x_c = 0.25, y_c =0.25$ on a domain which varies from $x\in[0,1]$ and $y\in [0,1]$. A total of $100 \times 100$ control volumes are considered in $x$ and $y$ directions and the solution is evolved using a time accurate scheme till $t=0.36$.
\begin{figure*}[!p]
\centering
\begin{multicols}{2}
{RICCA} \par {MOVERS+}
\end{multicols}
\begin{multicols}{2}
    \includegraphics[scale=0.35]{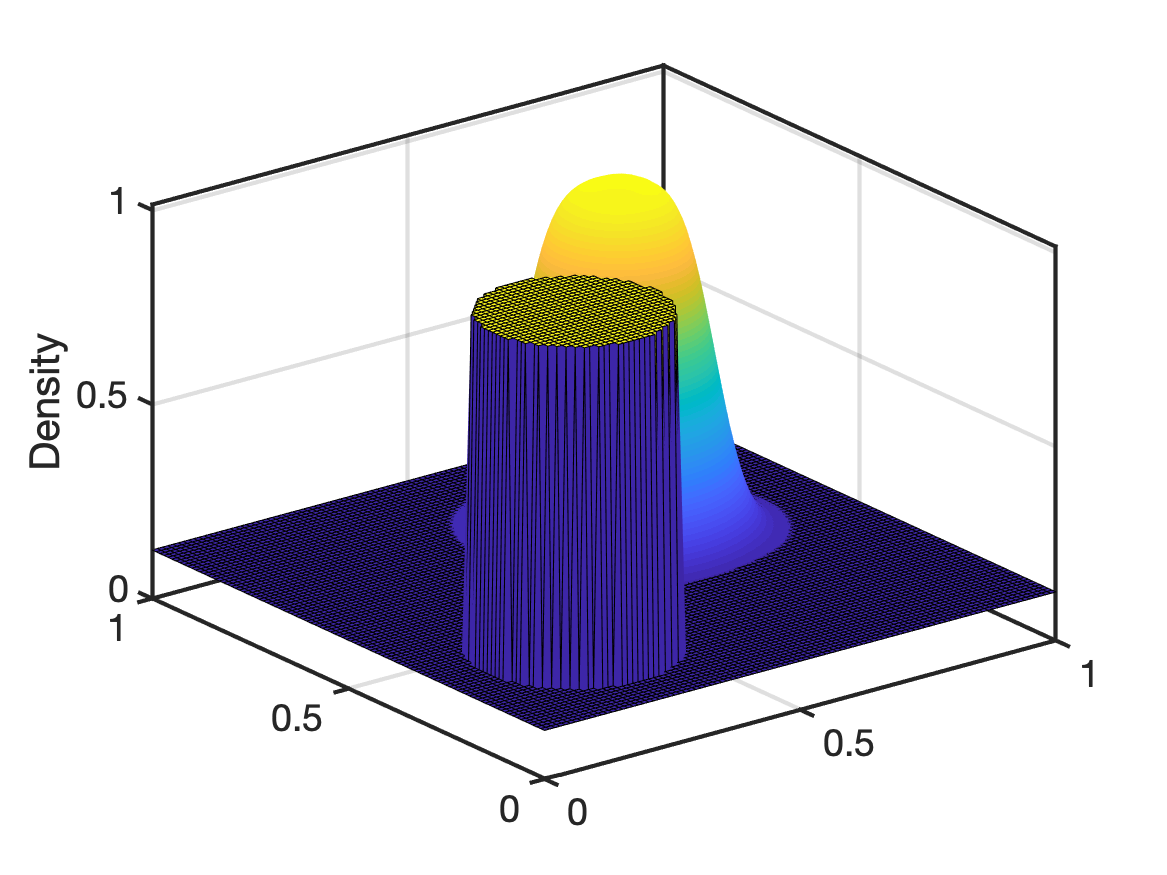}\par 
  \label{fig:RICCA_Interface_Only_PGEOS_100}
    \includegraphics[scale=0.35]{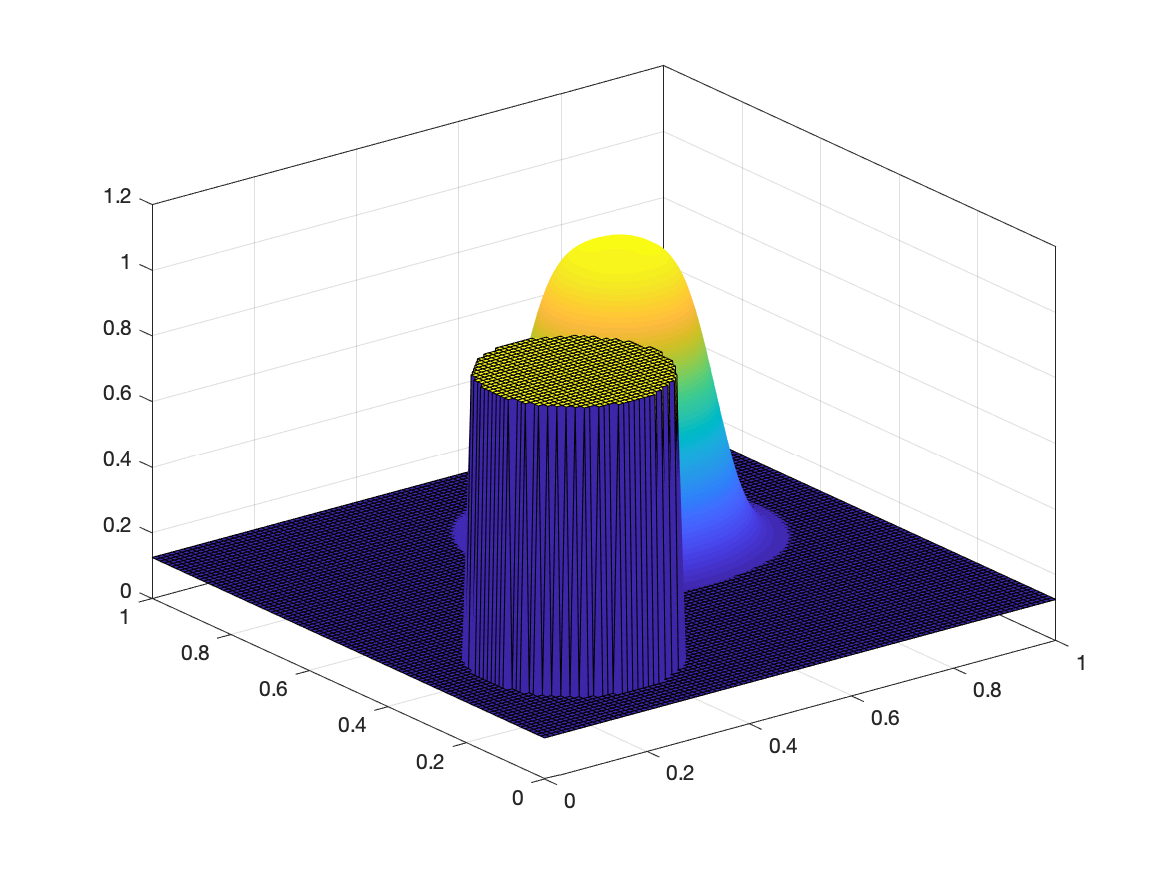}\par 
  \label{fig:Interface_Only_Pressure_PGEOS_100}
    \end{multicols}
    \caption{Surface of density at $t=0$ and $t =0 .36$}
    \begin{multicols}{2}
    \includegraphics[scale=0.35]{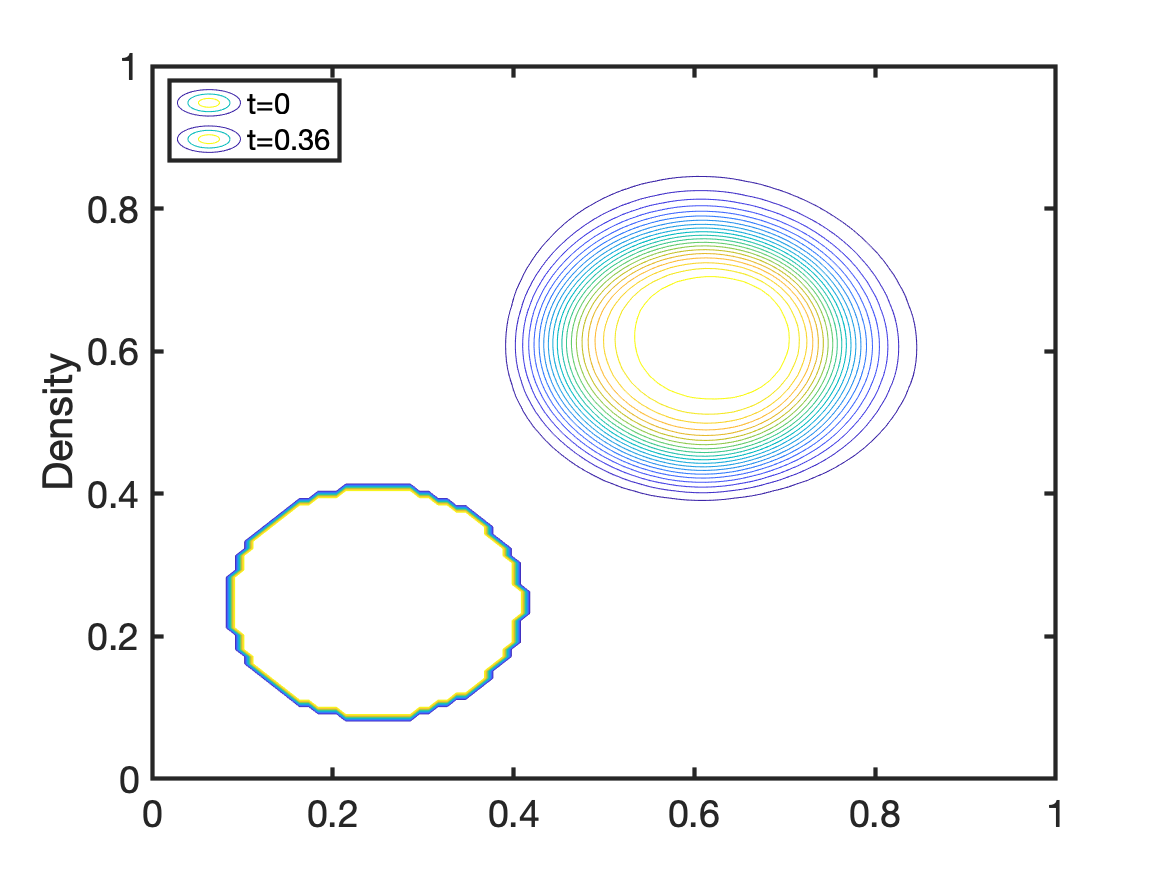}\par
  \label{Interface_Only_EOS_100}
    \includegraphics[scale=0.35]{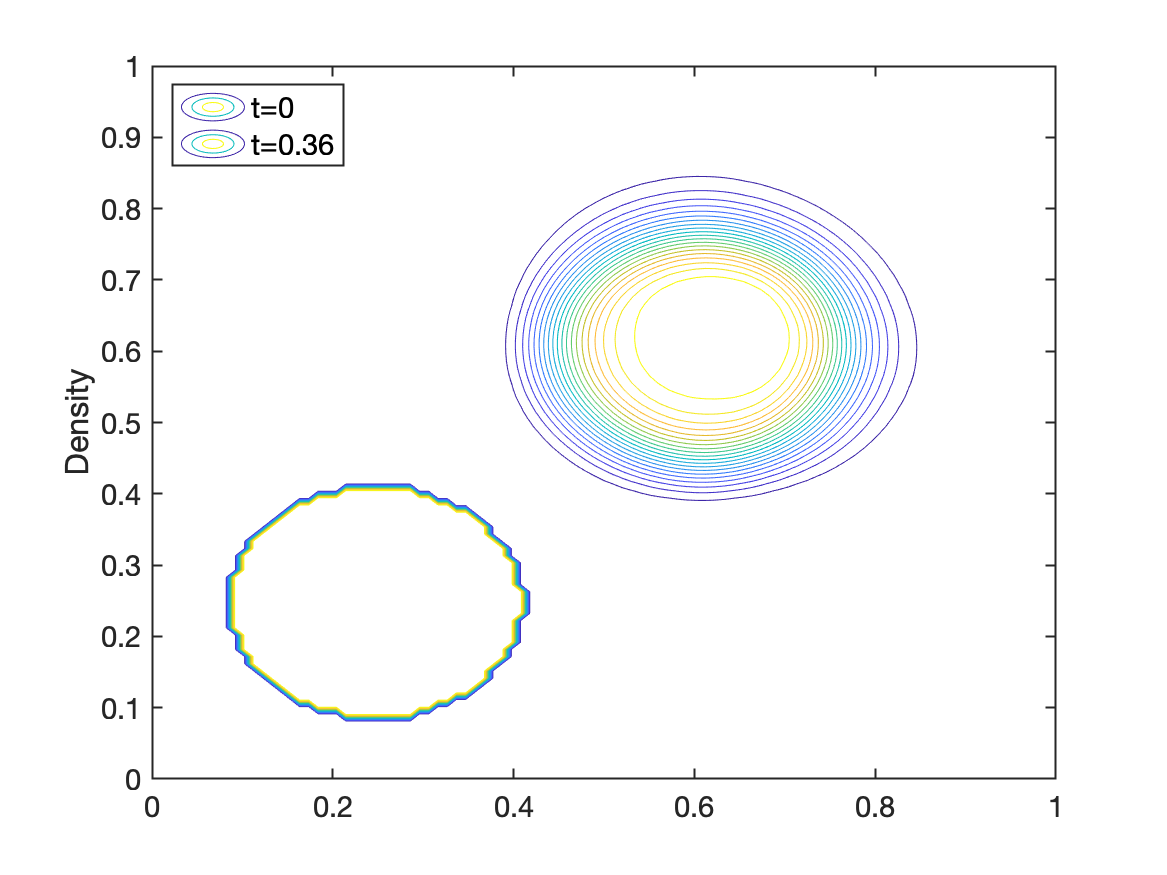}\par
  \label{Interface_Only_Density_PGEOS_100}
\end{multicols}
\caption{Density contours at $t=0$ and $t= 0.36$}
\begin{multicols}{2}
    \includegraphics[scale=0.35]{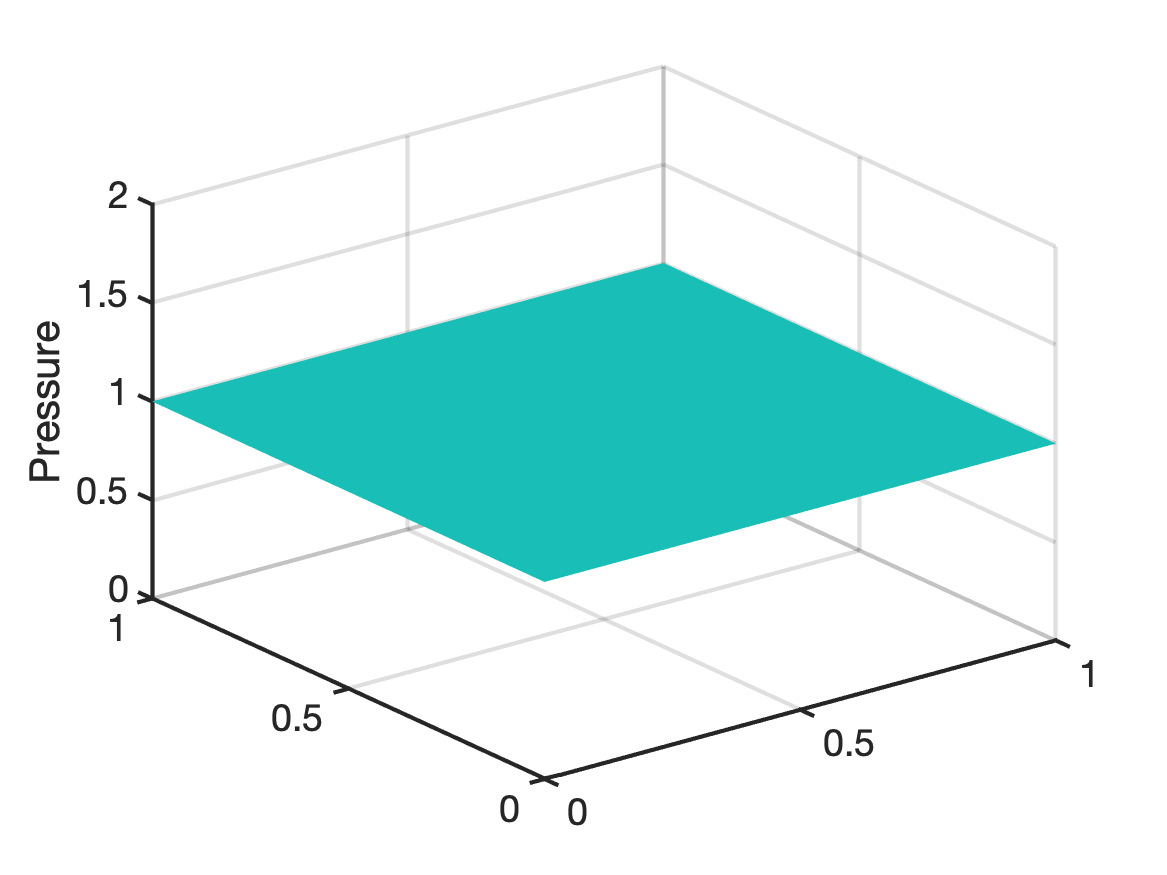}\par
  \label{Interface_Only_EOS_100_RICCA}
    \includegraphics[scale=0.35]{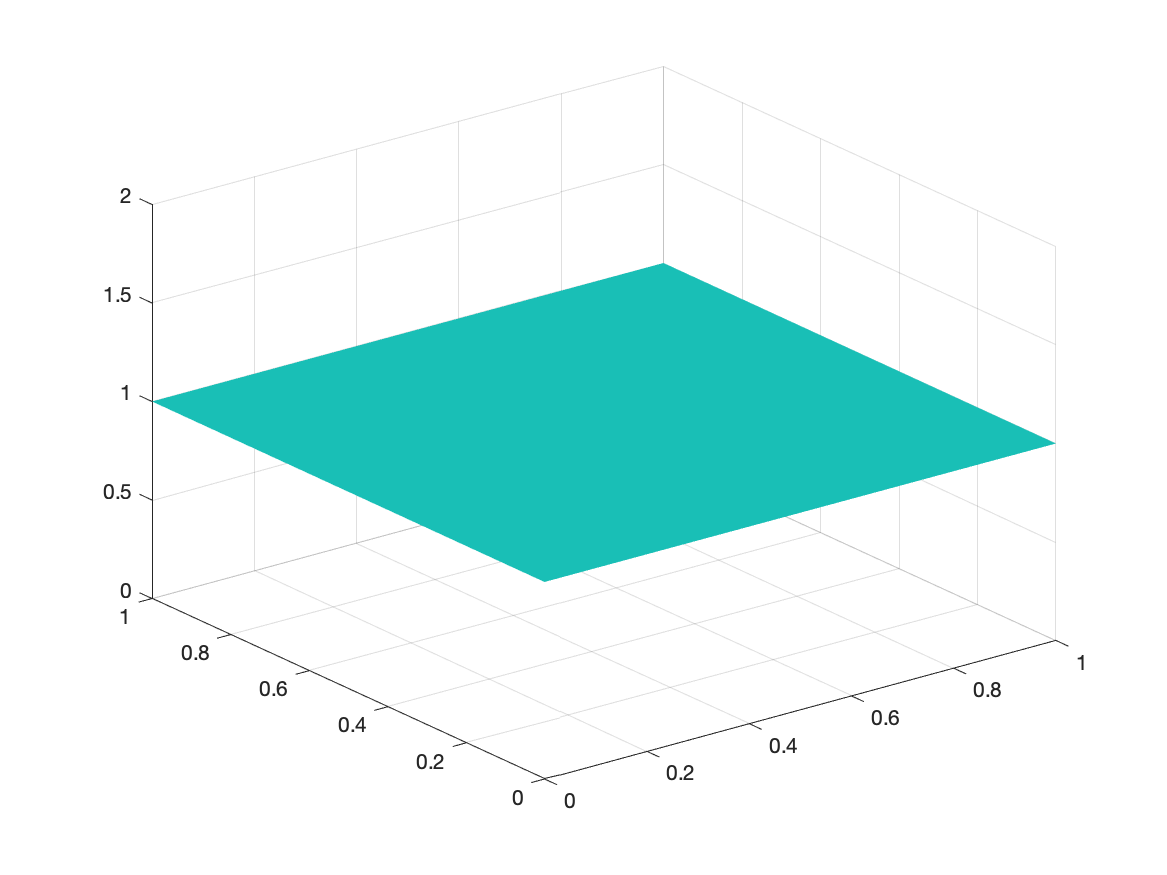}\par
  \label{Interface_Only_Pressure_PGEOS_100}
\end{multicols}
  \caption{Pressure distribution in the domain}
    \caption{Contour and surface view of density and pressure of interface only at $t = 0.36$ on $100 \times 100$ Grid}
    \label{Moving_Interface_100}
\end{figure*}

\begin{figure*}[!p]
\centering
\begin{multicols}{2}
{RICCA} \par {MOVERS+}
\end{multicols}
\begin{multicols}{2}
    \includegraphics[scale=0.35]{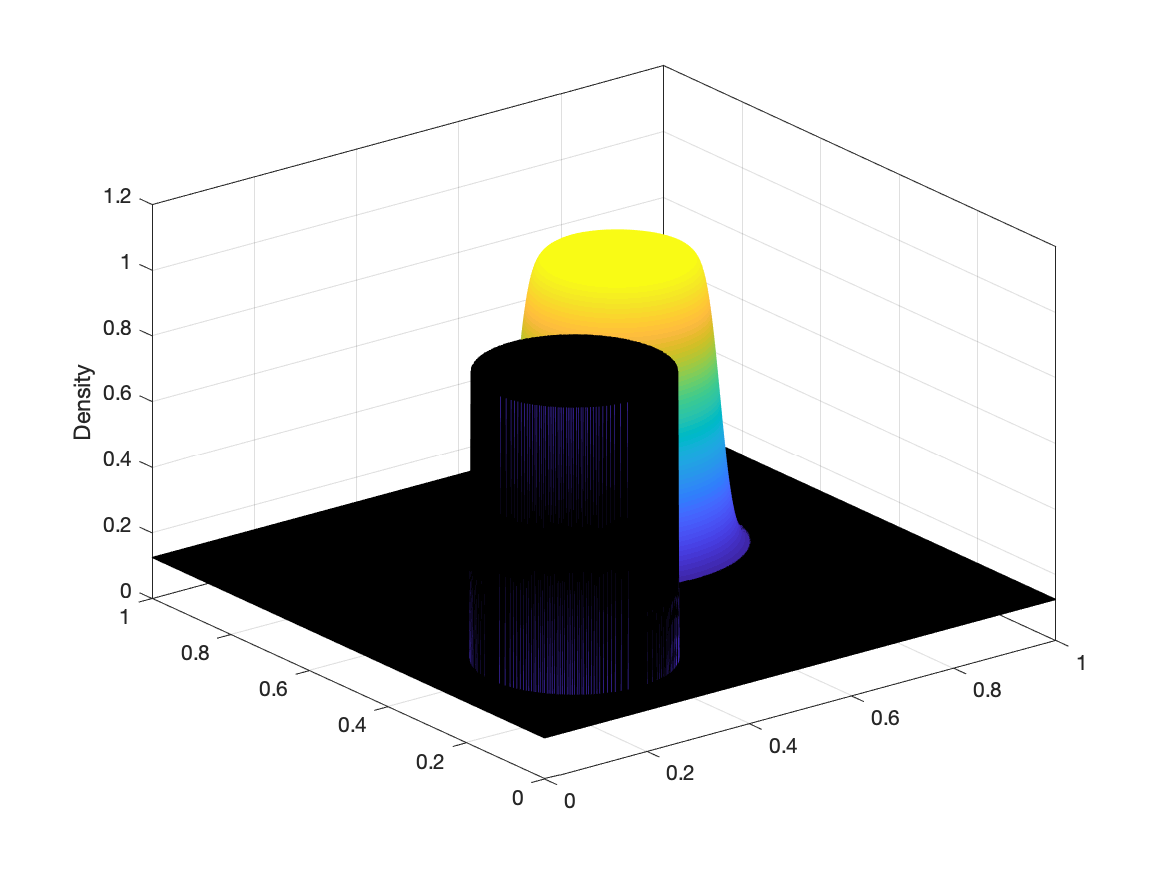}\par 
  \label{fig:RICCA_Interface_Only_PGEOS_500}
    \includegraphics[scale=0.35]{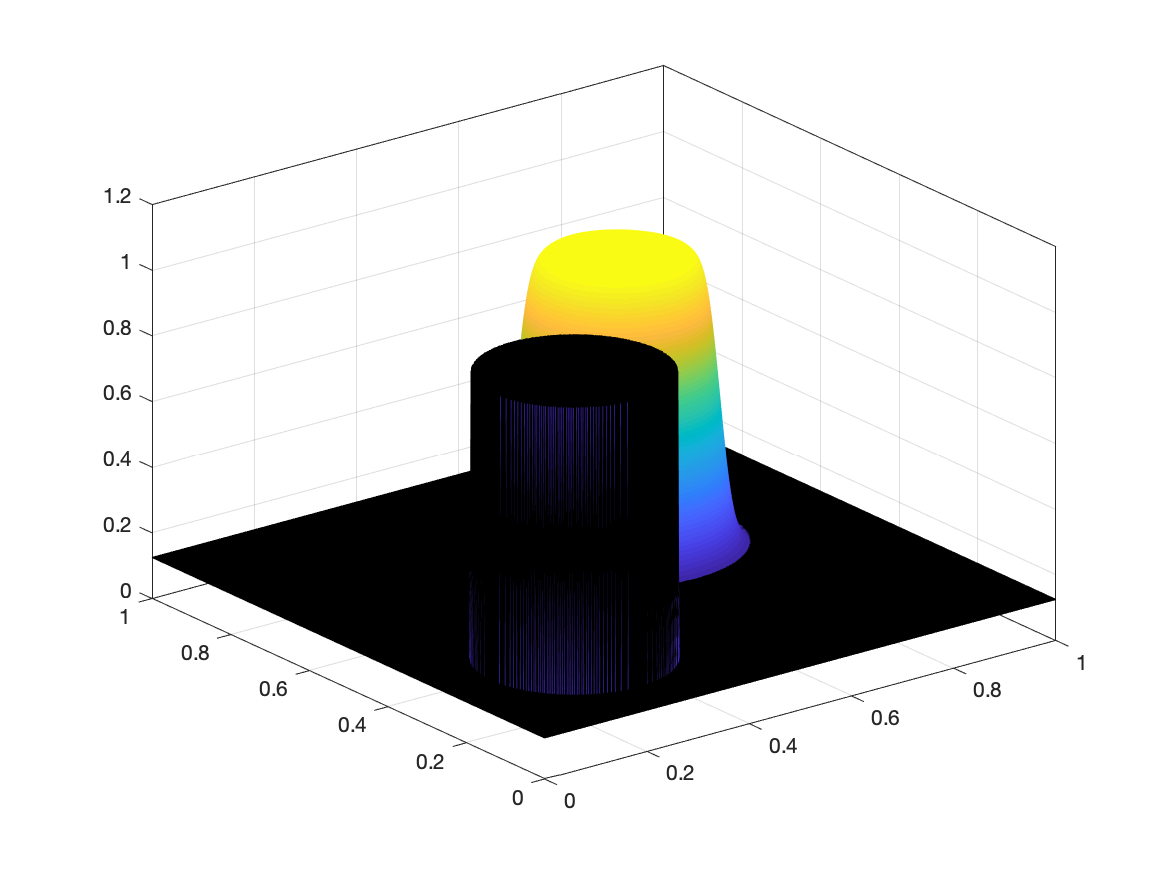}\par 
  \label{fig:NWSC_Interface_Only_Pressure_PGEOS_500}
    \end{multicols}
    \caption{Surface of density at $t=0$ and $t =0 .36$}
    \begin{multicols}{2}
    \includegraphics[scale=0.35]{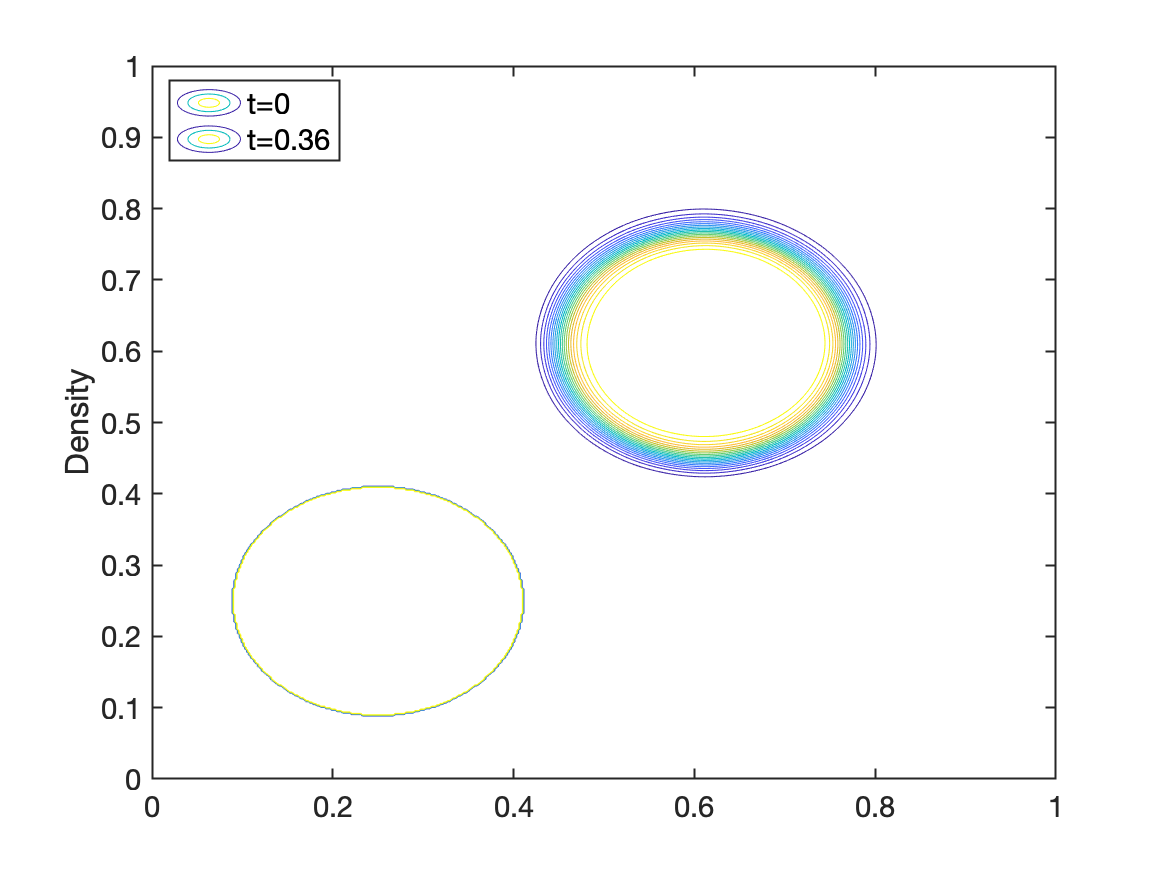}\par
  \label{RICCA_Interface_Only_Density_500}
    \includegraphics[scale=0.35]{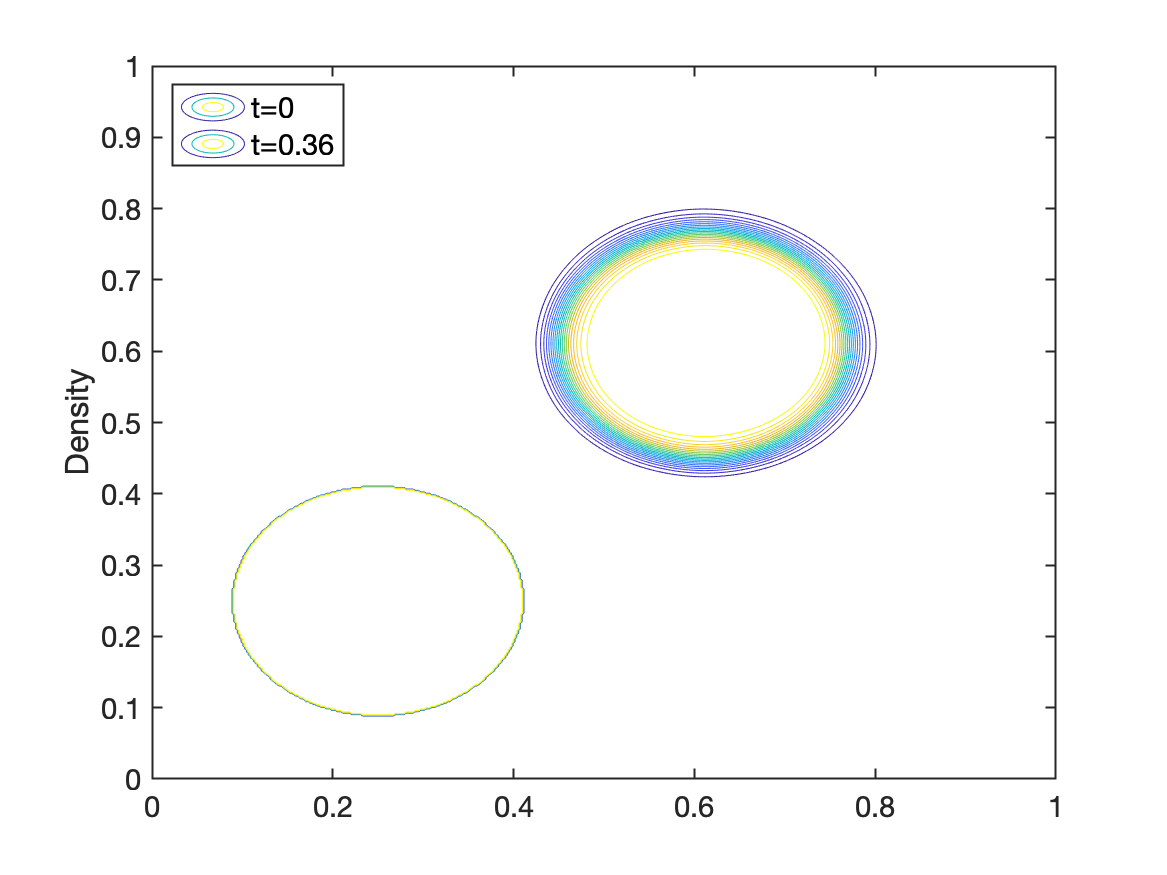}\par
  \label{NWSC_Interface_Only_Density_500}
\end{multicols}
\caption{Density contours at $t=0$ and $t= 0.36$}
\begin{multicols}{2}
    \includegraphics[scale=0.35]{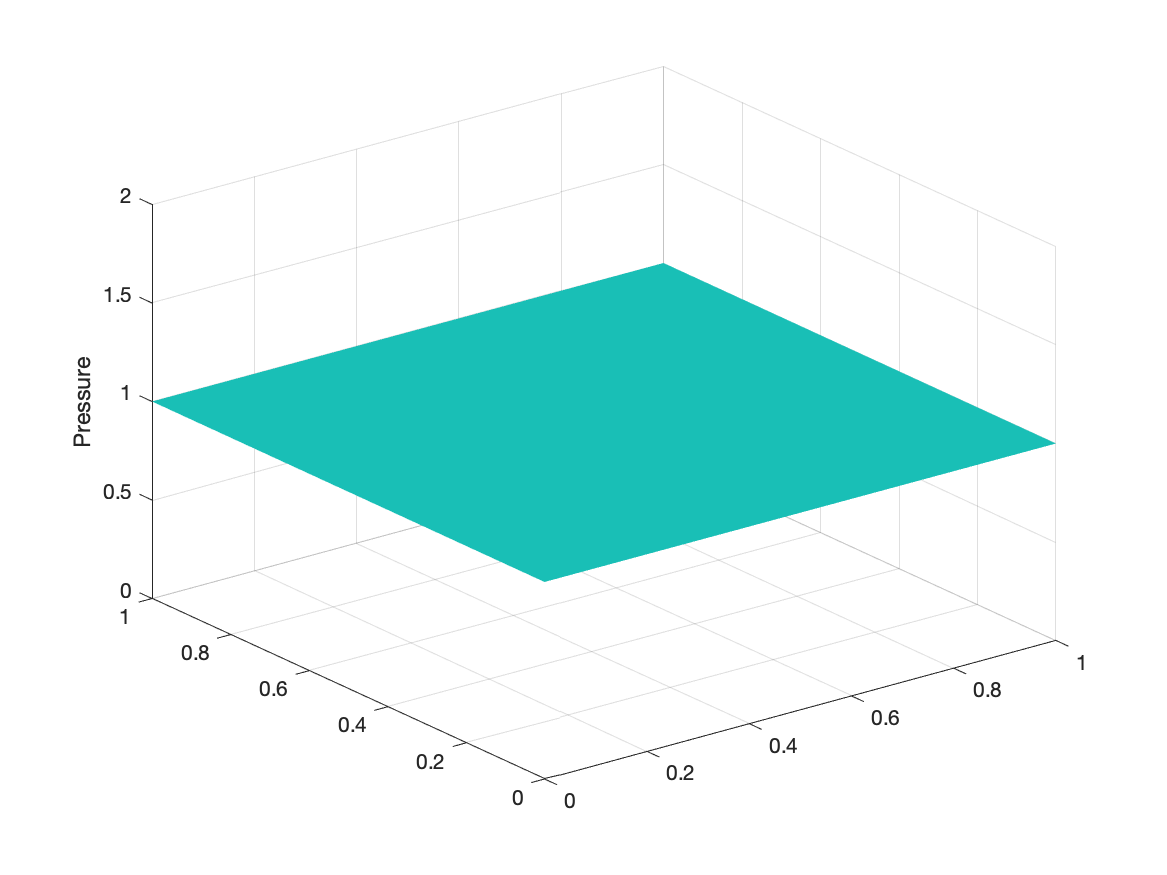}\par
  \label{RICCA_Interface_Only_EOS_Pressure_500}
    \includegraphics[scale=0.35]{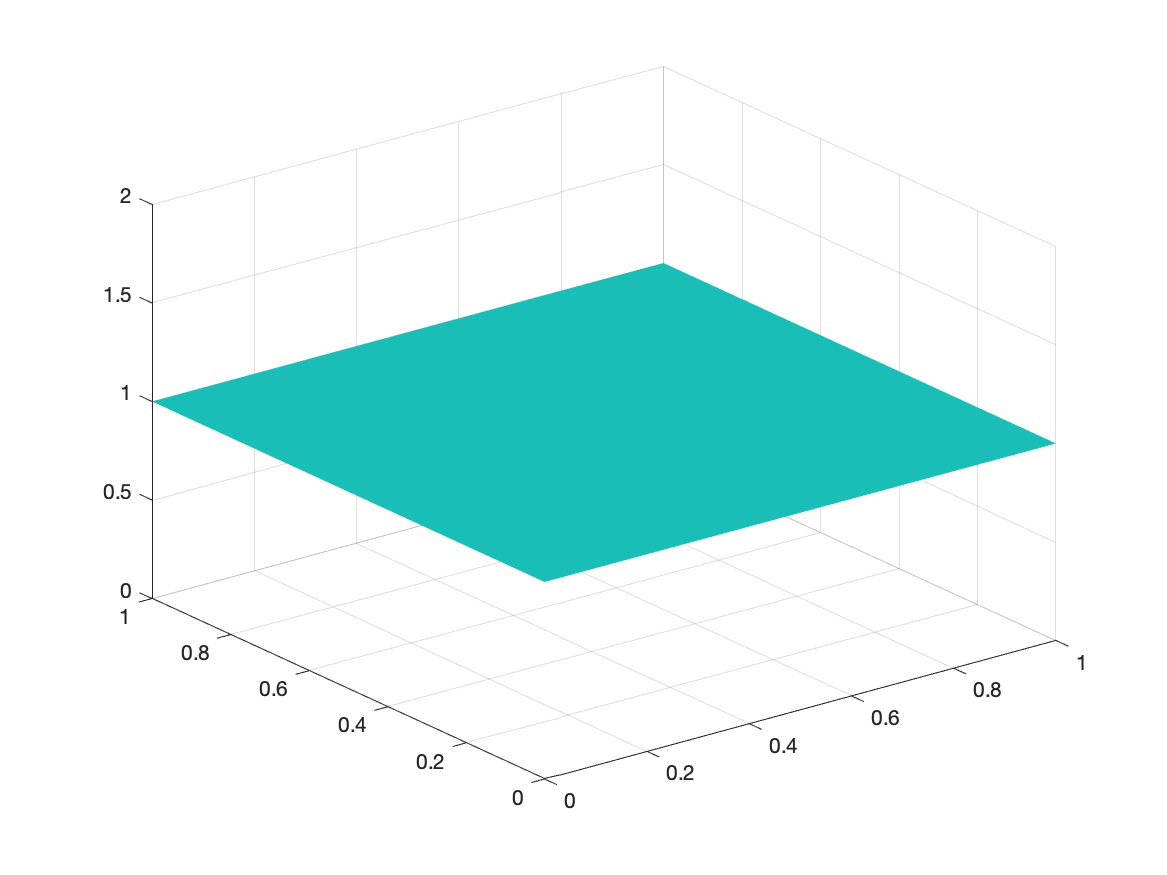}\par
  \label{NWSC_Interface_Only_Pressure_PGEOS_500}
\end{multicols}
  \caption{Pressure distribution in the domain}
    \caption{Contour and surface view of density and pressure of interface only at $t = 0.36$ on $500 \times 500$ Grid}
    \label{Moving_Interface_500}
\end{figure*}

\begin{figure}[h!]
\begin{center}
\includegraphics[scale=0.3]{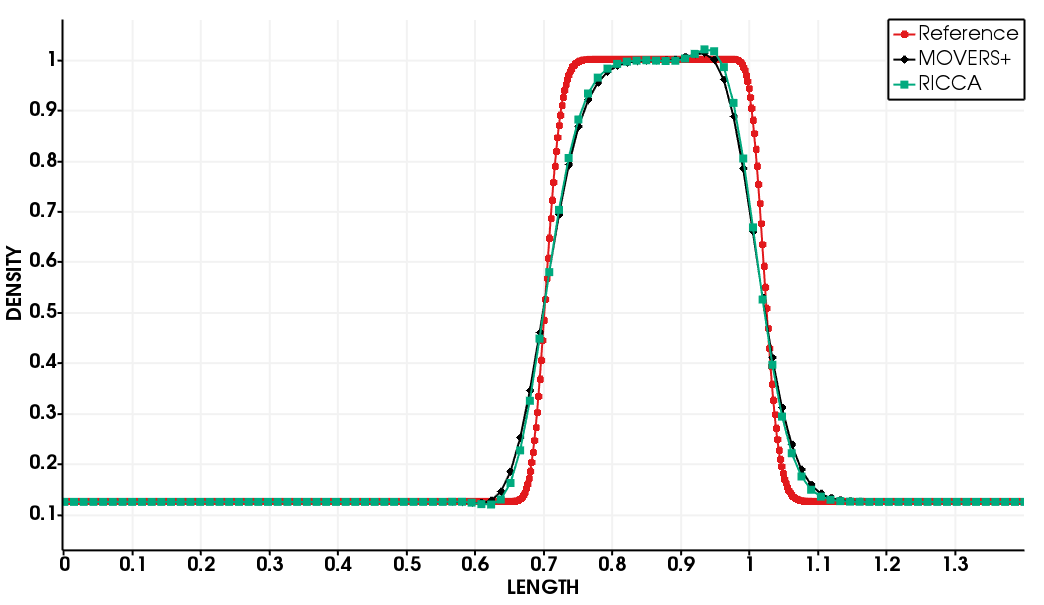}
\caption{Line plot of density across the interface}
\end{center}
\end{figure}
\begin{figure}[h!]
\begin{center}
\includegraphics[scale=0.3]{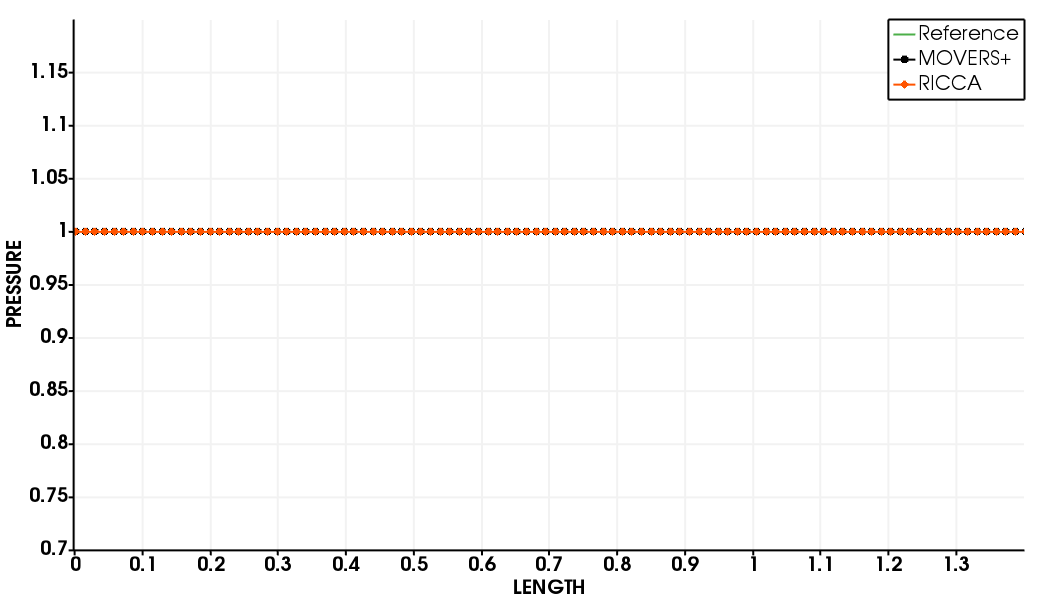}
\caption{Line plot of pressure across the interface}
\end{center}
\caption{2D Line plots of density and pressure across the interface using RICCA Scheme on $500\times 500$ and $100 \times 100$ Grid}
\label{fig:Line Plots}
\end{figure}
Figure (\ref{Moving_Interface_100}) refers to the moving bubble at $t=0$ and at $t=0.36$ using RICCA and MOVERS+ on a $100 \times 100$ grid and Figure (\ref{Moving_Interface_500}) refers to results obtained on a $500 \times 500$ grid. It can be observed that the position of the bubble is captured accurately and the pressure does not have any oscillations. Further the pressure and density plots across the bubble are shown in the figure (\ref{fig:Line Plots}) where in it can be observed that there are no pressure oscillations in the results generated by both RICCA and MOVERS+.

\subsubsection{Bubble explosion problem}
\begin{figure}[p!]
\begin{center}
    \includegraphics[scale = 0.75]{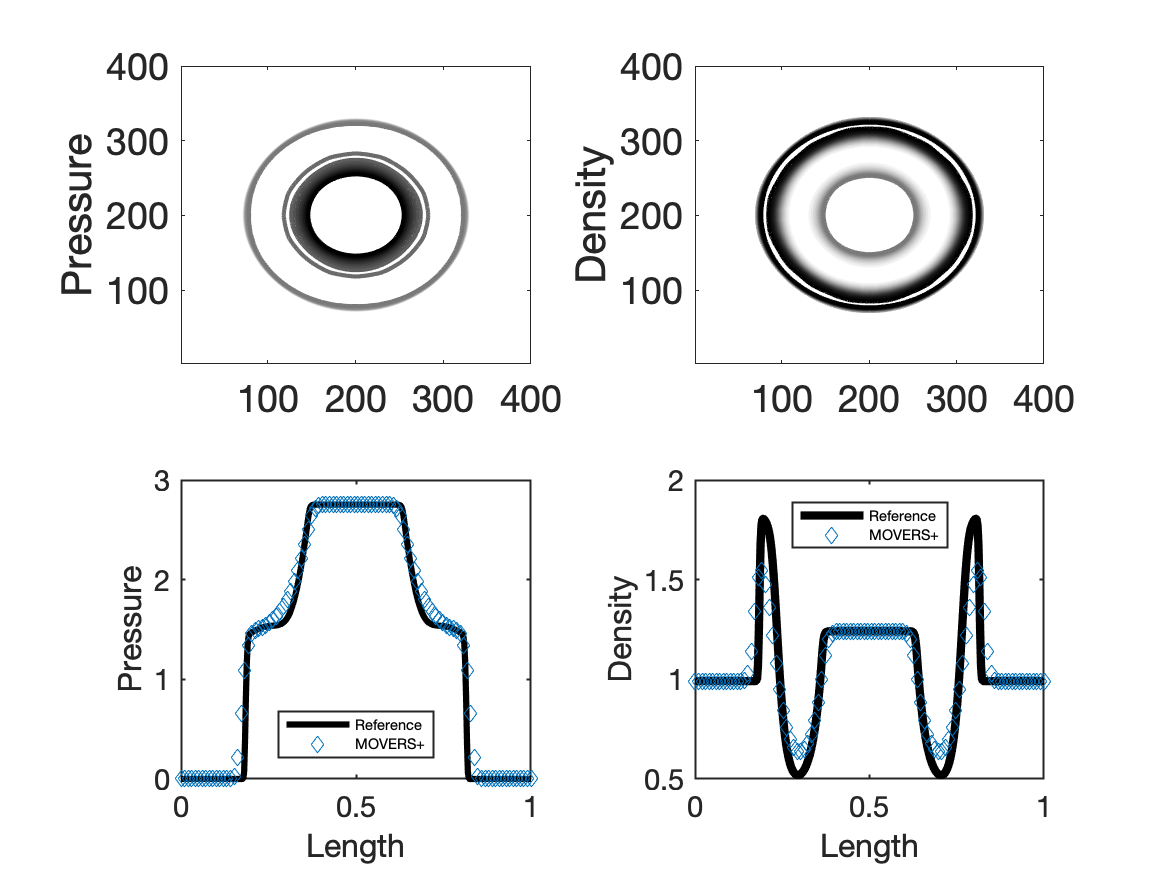}\par
    \caption{2D Line and contour plots of density and ressure using MOVERS+}
    \label{fig:2DPlots_BE_MOVERSNWSC}
    \includegraphics[scale=0.75]{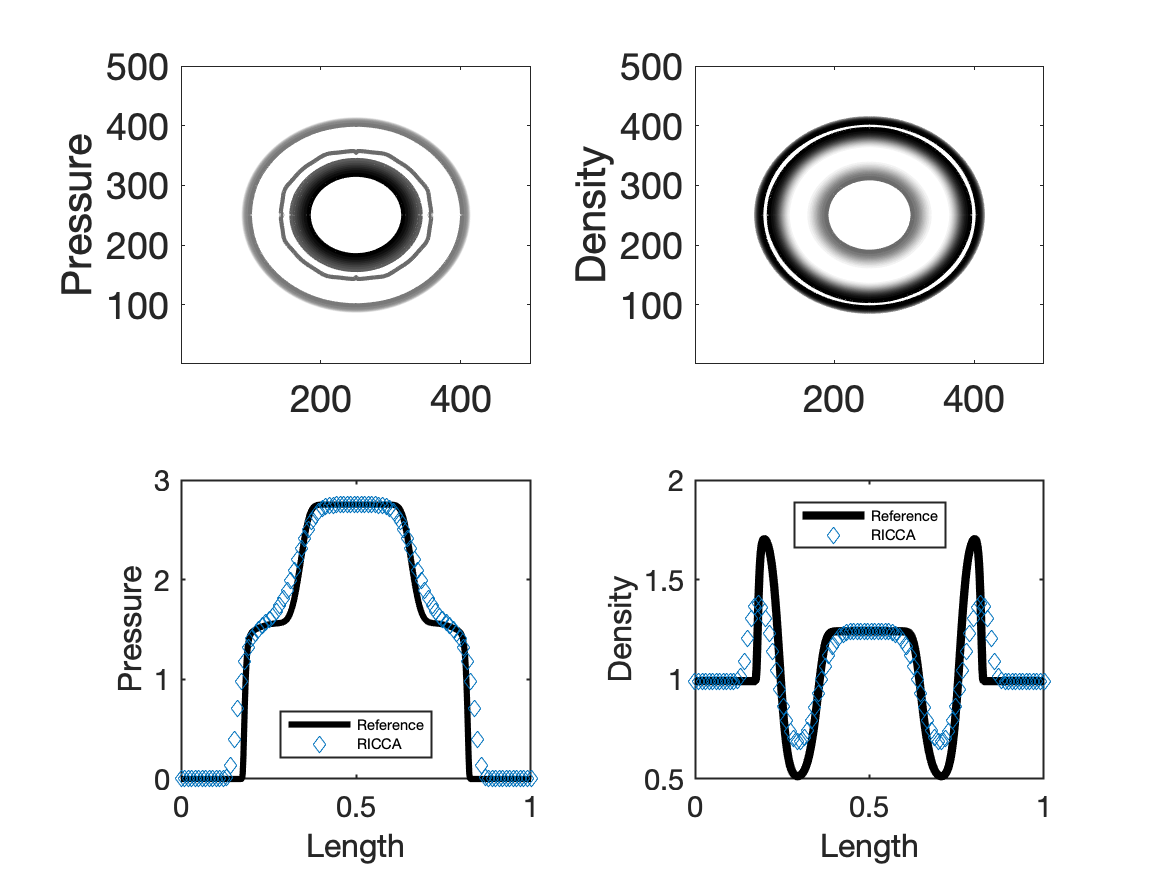}\par
    \caption{2D Line and contour plots of density and pressure using RICCA}
  \label{Bubble_Explosion_Contour_PGEOS_RICCA}
  \end{center}
\end{figure}

The second test case considered is a radially symmetric problem. It consists of a circular bubble present initially at rest in water and suddenly explodes due to high pressure of the water. The bubble is placed at $(x,y) = (0.5,0.5)$ and has a radius of $r_o = 0.2$. The fluid inside the bubble has the following initial conditions $(\rho, p, \gamma, p_{\infty}) = (1.241,2.753,1.4,0)$ and the surrounding water has the following properties $(\rho, p, \gamma, p_{\infty}) = (0.991,3.059\times 10^{-4},5.5,1.505)$. Numerical simulations are carried out using RICCA and MOVERS+. Pressure and density contour plots are shown in figure (\ref{fig:2DPlots_BE_MOVERSNWSC}) and the variation of pressure and density at $y=0.5$ are shown in figure (\ref{Bubble_Explosion_Contour_PGEOS_RICCA}).

\begin{figure*}[!p]
\centering
\begin{multicols}{2}
{MOVERS+} \par {RICCA}
\end{multicols}
\begin{multicols}{2}
    \includegraphics[scale=0.5]{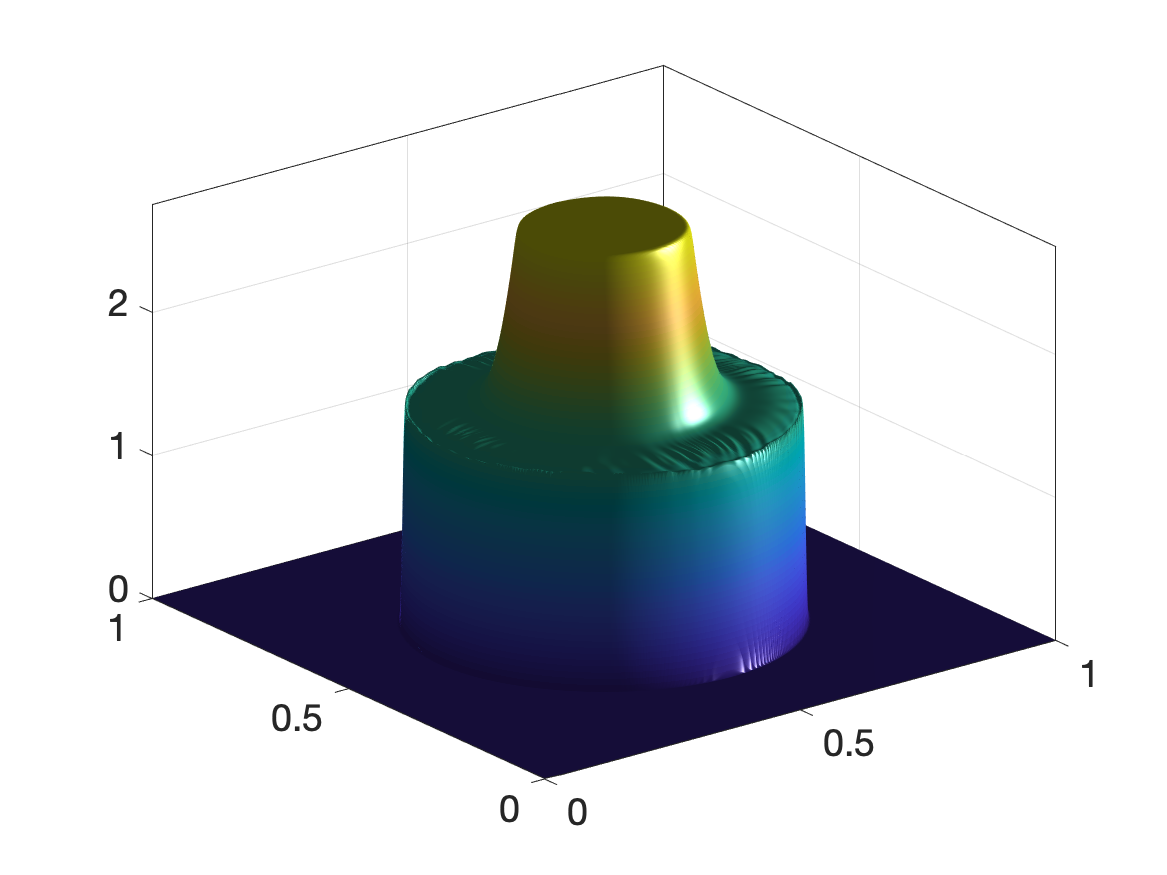}\par 
  \label{fig:RICCA_Bubble_Explosion_PGEOS}
    \includegraphics[scale=0.5]{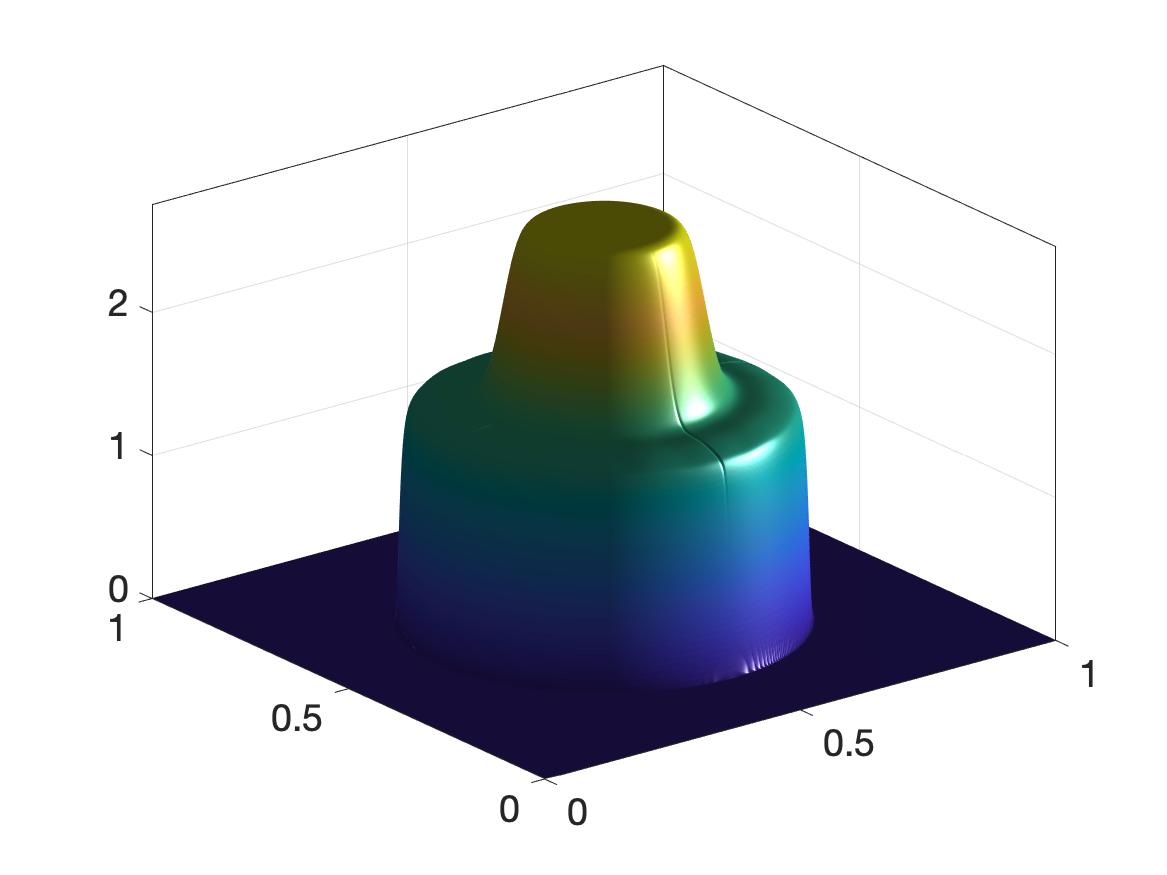}\par 
  \label{fig:Bubble_Explosion_Pressure_PGEOS}
    \end{multicols}
    \caption{Surface pressure of bubble explosion}
\begin{multicols}{2}
    \includegraphics[scale=0.5]{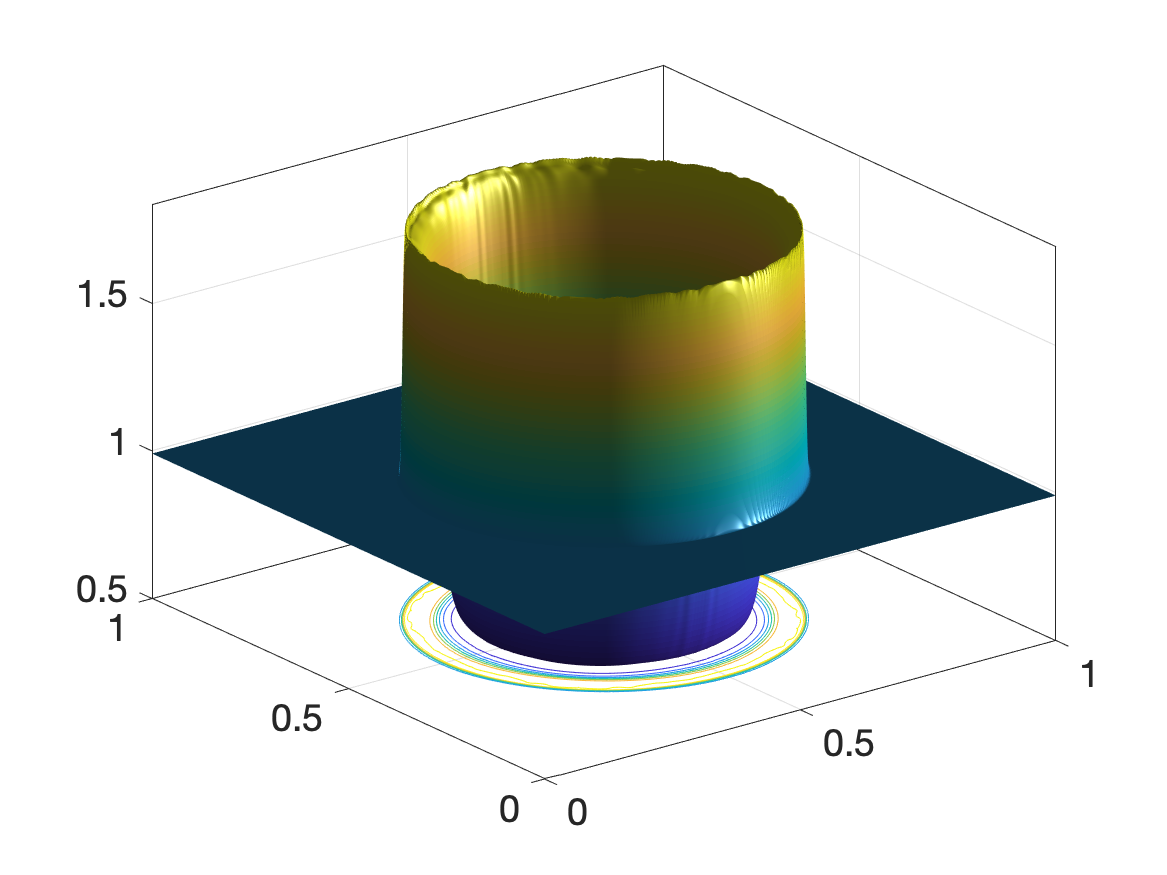}\par
  \label{Bubble_Explosion_EOS}
    \includegraphics[scale=0.5]{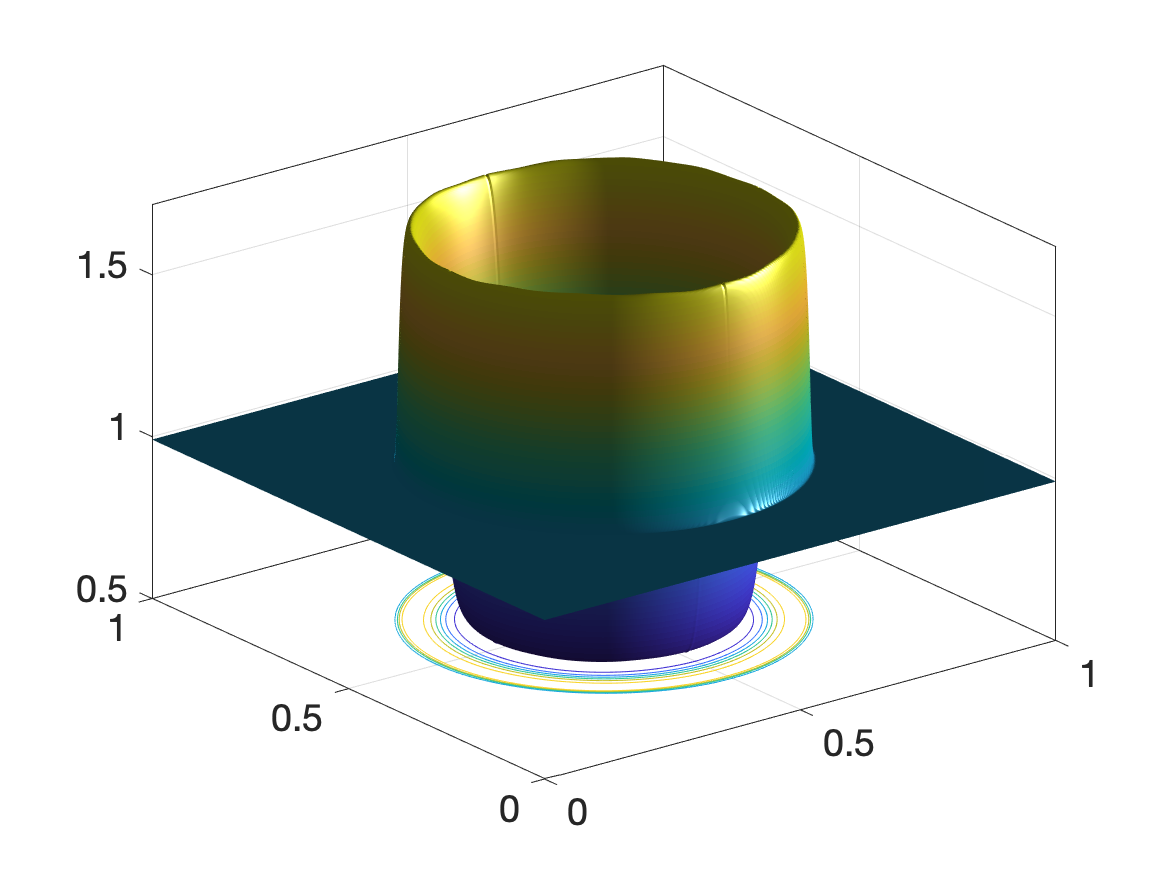}\par
  \label{Bubble_Explosion_Density_PGEOS}
\end{multicols}
    \caption{Surface density of bubble explosion}
    \caption{Contour and surface view of density and pressure of bubble explosion under water at $t = 0.058$ using MOVERS+ and RICCA on $500 \times 500$ Grid}
\end{figure*}

\section{Conclusions}\label{sec5}
Numerical simulations of Euler equations in 1D and 2D have been carried out using MOVERS (both scalar dissipation and vector dissipation), MOVERS+ and RICCA, for multicomponent gases with perfect gas EOS and stiffened gas EOS using mass fraction approach and $\gamma$-based approach in conservative form. It can be concluded that 
\begin{enumerate}
\item both the numerical schemes RICCA and MOVERS+ can be extended to multicomponent gases with different EOS without any modifications,
\item both the numerical schemes preserve the mass fraction positivity and the pressure positivity in the conservative approach when used in mass fraction based model.
\item pressure oscillations are observed in the finite volume framework for interface only problem and when $\gamma$-based model is adopted no pressure oscillations are observed.
\item both these schemes can be easily extended to any number of components.
\item both the schemes can handle large jumps in $\gamma$ without any modifications.
\end{enumerate}

\nocite{*}
\bibliography{wileyNJD-AMA}%



\end{document}